\documentclass[3p]{elsarticle}

\usepackage[colorlinks=true]{hyperref}
\usepackage[english]{babel}
\usepackage[utf8]{inputenc}

\usepackage{amsmath}
\usepackage{mathrsfs}
\usepackage{algorithm,algpseudocode}
\usepackage[figbotcap]{subfigure}
\usepackage{enumitem}
\usepackage{MnSymbol}%
\usepackage{wasysym}%
\usepackage{capt-of}
\usepackage[toc,page]{appendix}
\usepackage{accents}        
\usepackage[normalem]{ulem} 
\usepackage{tikz}            
\usepackage{empheq}         
\usepackage{todonotes}
\usepackage{placeins}
\usepackage{verbatim}
\usepackage{etoolbox}

\makeatletter
\newlength\appendixwidth
\preto\appendix{\addtocontents{toc}{\protect\patchl@Section}}
\newcommand{\patchl@Section}{%
  \settowidth{\appendixwidth}{\textbf{Appendix }}%
  \addtolength{\appendixwidth}{1.5em}%
  \patchcmd{\l@Section}{1.5em}{\appendixwidth}{}{\ddt}%
}
\makeatother

\newcommand\acclrvec[1]{\accentset{\,\leftrightarrow}{#1}}	
\newcommand\norm[1]{\left\lVert#1\right\rVert}
\newcommand{\tauvec}{\pmb{\tau}}

\newcommand{\epsvec}{\pmb{\epsilon}}

\def\Res{\pmb{\mathfrak{R}}}
\def\F{\pmb{\mathfrak{F}}}
\newcommand{\numflux}[1]{\hat{\mathbf{#1}} } 
\newcommand{\blocktensor}[1]{\acclrvec{{\mathbf #1}}\ }		
\newcommand{\Nabla} {\vec{\nabla}}
\def\d{\mathrm{d}}
\newcommand{\mat}[1]{\uline{\mathbf{#1}}}
\newcommand{\bigderiv}[2]{ \frac{\d {#1}}{\d {#2} } }

\newcommand{\NDOF} {\mathrm{NDOF}}
\def\Re {\mathrm{Re}_{\infty}}%
\def\Ma {\mathrm{Ma}_{\infty}}
\newcommand{\gstate}[1]{\boldsymbol{#1}}

\journal{Journal of Applied Mathematics and Computation}

\biboptions{sort&compress}

\begin{document}

\tikzstyle{endbox} = [rectangle, rounded corners, minimum width=3cm, minimum height=1cm,text centered, draw=black, fill=red!30]
\tikzstyle{startbox} = [rectangle, rounded corners, minimum width=3cm, minimum height=1cm,text centered, draw=black, fill=green!30]
\tikzstyle{process} = [rectangle, minimum width=3cm, minimum height=1cm, text centered, draw=white, fill=white]
\tikzstyle{decision} = [rectangle, minimum width=3cm, minimum height=1cm, text centered, draw=white, fill=white]
\tikzstyle{point} = [coordinate, node distance=0.5cm]
\tikzstyle{arrow} = [thick,->,>=stealth]

\hypersetup{
urlcolor=black
}

\begin{frontmatter}
\title{Truncation Error-Based Anisotropic $p$-Adaptation for Unsteady Flows for High-Order Discontinuous Galerkin Methods}

\author[COLOGNE]{Andrés M. Rueda-Ramírez \corref{mycorrespondingauthor}}
\cortext[mycorrespondingauthor]{Corresponding authors:}
\ead{aruedara@uni-koeln.de}

\author[DMAE]{Gerasimos Ntoukas}

\author[DMAE,CCS]{Gonzalo Rubio \corref{mycorrespondingauthor}}
\ead{g.rubio@upm.es}

\author[DMAE,CCS]{Eusebio Valero}

\author[DMAE,CCS]{Esteban Ferrer}

\address[COLOGNE]{Department of Mathematics and Computer Science, University of Cologne, 50931 Cologne, Germany}

\address[DMAE]{ETSIAE-UPM (School of Aeronautics), Universidad Politécnica de Madrid, Plaza de Cardenal Cisneros 3, 28040 Madrid, Spain}

\address[CCS]{Center for Computational Simulation, Universidad Politécnica de Madrid, Campus de Montegancedo, Boadilla del Monte, 28660 Madrid, Spain}

\begin{abstract}
In this work, we extend the $\tau$-estimation method to unsteady problems and use it to adapt the polynomial degree for high-order discontinuous Galerkin simulations of unsteady flows. The adaptation is local and anisotropic and allows capturing relevant unsteady flow features while enhancing the accuracy of time evolving functionals (e.g., lift, drag).
To achieve an efficient and unsteady truncation error-based $p$-adaptation scheme, we first revisit the definition of the truncation error, studying the effect of the treatment of the mass matrix arising from the temporal term. 
Secondly, we extend the $\tau$-estimation strategy to unsteady problems. 
Finally, we present and compare two adaptation strategies for unsteady problems: the dynamic and static $p$-adaptation methods. In the first one (dynamic) the error is measured periodically during a simulation and the polynomial degree is adapted immediately after every estimation procedure. In the second one (static) the error is also measured periodically, but only one $p$-adaptation process is performed after several estimation stages, using a combination of the periodic error measures. The static $p$-adaptation strategy is suitable for time-periodic flows, while the dynamic one can be generalized to any flow evolution.    

We consider two test cases to evaluate the efficiency of the proposed $p$-adaptation strategies.
The first one considers the compressible Euler equations to simulate the advection of a density pulse. The second one solves the compressible Navier-Stokes equations to simulate the flow around a cylinder at Re=100.    
The local and anisotropic adaptation enables significant reductions in the number of degrees of freedom with respect to uniform refinement, leading to speed-ups of up to $\times4.5$ for the Euler test case and $\times2.2$ for the Navier-Stokes test case.

\end{abstract}

\begin{keyword}
High-order discontinuous Galerkin, Anisotropic $p$-adaptation, Unsteady $p$-adaptation, Compressible flows.
\end{keyword}

\end{frontmatter}


\section{Introduction}

High-order DG methods are expected to be the engine of the next generation of CFD codes \cite{wagner2022performance}. Commonly, high-order DG schemes are formulated as multi-domain spectral methods. As a result, besides the increased accuracy of spectral (high-order) methods, they also provide a compact stencil and therefore a local character, a feature that makes them highly parallelizable and flexible for complex 3D geometries \cite{Cockburn2000,Wang2013High}. Moreover, DG methods can handle non-conforming meshes with hanging nodes and/or different polynomial degrees efficiently \cite{Ferrer2012a,FERRER2017754,ferrer2012high}, which makes them well suited for mesh adaptation strategies.

Among the different DG formulations, the discontinuous Galerkin spectral element method (DGSEM) \cite{Black1999,kopriva2009implementing} is a nodal (collocation) version of the DG method that uses tensor-product Lagrange basis functions (traditionally in hexahedra) and stores the data at the Gauss or Gauss-Lobatto nodes of the quadrature rule. The use of a quadrature rule with the same number of nodes as the approximate solution equips the DGSEM with a diagonal mass matrix and very cheap-to-compute operators. In fact, the computational cost of the DGSEM has been estimated to be a factor of four smaller than other conventional DG methods \cite{Beck2016}.
In addition, since the DGSEM uses tensor-product bases, it can handle $p$-anisotropic discretizations efficiently \cite{Kompenhans2016,RuedaRamirez2019b}, i.e. discretizations that have different polynomial degrees in each coordinate direction. For all those properties, the DGSEM has been used in a wide range of applications, including the simulation of incompressible Navier-Stokes, compressible Navier-Stokes, Cahn Hilliard equation, or multi-phase flows (see \cite{ferrer2022horses3d} and references therein).

In their famous review paper, Wang et al. \cite{Wang2013High} point out that one of the challenges the high-order community must address to impact the design process and replace traditional low-order codes is the development of efficient mesh adaptation strategies. The idea behind these strategies is to reduce the number of degrees of freedom (DOFs) while maintaining high accuracy, which translates into shorter computational times and reduced storage requirements. 
Local adaptation can be performed by subdividing or merging elements ($h$-adaptation), by enriching or reducing the polynomial degree in certain elements ($p$-adaptation), by relocating the position of the nodes in a mesh ($r$-adaptation). For all these strategies it is of paramount importance to identify the flow regions that require refinement or coarsening with a local error estimation.

Adaptation strategies have been classified according to the type of error measure that is employed as feature-based adaptation, adjoint-based adaptation, and local error-based adaptation. A comparison of these three approaches was performed by Fraysse et al. \cite{Fraysse2012} for finite volume approximations and by Kompenhans et al. \cite{Kompenhans2016a} and Naddei et al. \cite{Naddei2018} for high-order DG methods.
The feature-based adaptation is the classical approach and uses easy-to-compute error measures that depend on the flow features. They rely on the assumption that high errors are expected where the flow is more difficult to resolve. Hence, refinement is predicted where high velocity, density or pressure gradients are identified \cite{Aftosmis1994,Persson2006}. For DG discretizations, an easy-to-compute feature-based adaptation criterion is the assessment of jumps across element interfaces \cite{Krivodonova2003,Krivodonova2004,Remacle2003}.
The main disadvantage of these methods is that there is no direct relation between the adaptation criterion and the numerical errors and thus the accuracy is not easily predictable. Additionally, the only way to solve steady-state problems is to adapt iteratively.

A second and more sophisticated approach is known as adjoint-based adaptation. In this approach, a functional target is defined (e.g. drag or lift in external flow aerodynamics) and the adjoint problem is solved to obtain a spatial distribution of the functional error, which is then used to adapt the mesh. This technique was originally developed for structural analysis using FEM by Babuška and Miller \cite{babuvska1984,babuvska1984a}, and has been used recently for adaptation strategies in DG methods \cite{Hartmann2006,Hartmann2002,Wang2009}.
The main drawback of this approach is the high computational cost to solve the adjoint problem and the storage requirements needed to save the error estimators, especially in unsteady flows. Moreover, only the error of the functional analyzed is guaranteed to be reduced, whereas the error of other functionals may deteriorate.

A computationally more efficient alternative is the local error-based adaptation, which is based on the assessment of any measurable (not feature-based) local error in all the cells of the domain \cite{Hartmann2002}. The local error-based adaptation methods are interesting since, in contrast to feature-based methods, they provide a way to predict and control the overall accuracy, and are computationally cheaper than adjoint-based schemes \cite{Kompenhans2016,Kompenhans2016a}.
For those reasons, local error-based adaptation strategies are retained in this work.
A large amount of effort has been invested in the development of reliable local error-based adaptation methods. Estimations of the local discretization error have been used by Mavriplis \cite{Mavriplis1989,Mavriplis1994} to develop $hp$-adaptation techniques for the spectral element method. 
Residual-based $p$-adaptation is also a local error-based adaptation method, which uses the residual to measure how accurate is the local approximation. This method was originally developed for Finite Elements (FE) and has been successfully used with DG methods \cite{Hartmann2006,Naddei2018}.
In the case of modal (hierarchical) DG methods, a possibility is to employ low cost error estimates that take advantage of the modal approximation to drive $p$-adaptation procedures, such as the Variational Multiscale (VMS) indicator by Kuru and De la Llave Plata \cite{Kuru2016}, or the spectral decay indicator by Persson and Peraire \cite{Persson2006}. 

In this work, we favor truncation error estimators, another local error-based alternative to drive a mesh adaptation method. The truncation error is related to the discretization error through the Discretization Error Transport Equation \cite{Roy2010}, where it acts as a local source term. This relation makes it useful as an indicator for mesh adaptation methods \cite{Choudhary2013,Syrakos2012}, since refining the mesh where the truncation error is high reduces the discretization error in all the mesh \cite{Rubio2015}, with an additional advantage: truncation error estimation requires less computational effort than adjoint methods. Finally, it has been shown that controlling the truncation error targets the numerical accuracy of all functionals at once \cite{Kompenhans2016,laskowski2022functional}, ensuring that adapting a mesh using the truncation error leads necessarily to an error decrease in any other functional (e.g. lift or drag). 

The $\tau$-estimation method proposed by Brandt \cite{Brandt1984}, which estimates the local truncation error by injecting a fine grid solution into coarser meshes, has been used to perform local error-based mesh adaptation in low-order schemes 
\cite{berger1987adaptive,Fraysse2014,Fraysse2012,Fraysse2013,Syrakos2012,Syrakos2006}. Rubio et al. \cite{Rubio2013} extended the $\tau$-estimation approach to high-order methods using a continuous Chebyshev collocation method. Later, Rubio et al. \cite{Rubio2015} applied it to DGSEM discretizations. Kompenhans et al. \cite{Kompenhans2016} applied the $\tau$-estimation approach to perform steady-state $p$-adaptation using the Euler and Navier-Stokes equations, and showed that a reduction of the truncation error increases the numerical accuracy of all functionals at once. Furthermore, Kompenhans et al. \cite{Kompenhans2016a} also showed that truncation error-based adaptation can exhibit better performance than feature-based adaptation.
In contrast to most local error-based adaptation methods, where multiple error estimation and adaptation stages are needed in a steady-state solution, the $\tau$-estimation method generates a unique prediction of what polynomial degree is needed for a desired truncation error threshold. Therefore, in steady-state the adaptation strategy is to converge a high-order approximation (reference mesh) to a specified global residual and then to perform a single error estimation followed by a corresponding $p$-adaptation process. Besides, the truncation error is known to decay exponentially in smooth solutions \cite{Kompenhans2016,Rubio2015}. Therefore, if the estimation is good, it is possible to extrapolate the behavior and predict the polynomial degree needed for a desired error threshold \cite{Kompenhans2016,RuedaRamirez2019}.\\

Being a relatively recent technique, $p$-adaptation methods that use $\tau$-estimators have only been applied to steady-state solutions with high-order methods \cite{Kompenhans2016,RuedaRamirez2019b}. In this work, we propose a methodology to extend this technique to unsteady problems. The rest of this work is organized as follows. First,
in Section \ref{sec:Uadapt:NewTruncError}, we show that the truncation error can be formulated in several forms, depending on the choice of the continuous and discrete partial differential operators.
Herein, a thorough analysis of the formulation that is traditionally used in the DG community \cite{Rubio2015,Kompenhans2016,Kompenhans2016a,RuedaRamirez2019,
RuedaRamirez2019a} is presented, a new formulation for the truncation error is proposed, and a comparative analysis of both techniques is detailed.
Second, in Section \ref{sec:Uadapt:Estim} we provide two strategies to estimate the truncation error in unsteady problems, both of which are derived from the $\tau$-estimation method.
The first strategy is directly derived from the variational DG formulation while the second strategy uses a dual time-stepping pseudo-time discretization.
Third, in Section \ref{sec:Uadapt:pAdaptAlg}, two $p$-adaptation algorithms for unsteady problems are proposed, which use an unsteady $\tau$-estimation method: (i) a dynamic $p$-adaptation strategy, which performs several stages of estimation and $p$-adaptation throughout a simulation and second, and (ii) a static $p$-adaptation strategy, which performs several truncation error estimation stages, but only one $p$-adaptation stage.
Finally, the methods are applied to unsteady problems modeled by the compressible Euler and Navier-Stokes equations in Section \ref{sec:Uadapt:Results}, and a detailed analysis of their performance is presented.
The most important findings of this paper are summarized in Section \ref{sec:Uadapt:Conclustions}.

\section{Numerical Methods}

\subsection{The Discontinuous Galerkin Spectral Element Method}\label{sec:DGSEM}

We consider the approximation of systems of conservation laws, 
\begin{equation}\label{eq:NScons}
\partial_t \mathbf{q} + F(\mathbf{q}) = \mathbf{0}, \ \ \text{in } \Omega,
\end{equation}
subject to appropriate boundary conditions, where $\mathbf{q}$ is the state vector of conserved variables, and $F(\mathbf{q})=\Nabla \cdot \blocktensor{f}$ is the continuous partial differential operator, where $\blocktensor{f}$ is a flux block vector, which depends on $\mathbf{q}$. 

In an advection-diffusion conservation law, such as the Navier-Stokes equations, the flux vector can be written as
\begin{equation}
\blocktensor{f} = \blocktensor{f}^a(\mathbf{q}) - \blocktensor{f}^{\nu}(\mathbf{q}, \Nabla \mathbf{q}),
\end{equation}
where $\blocktensor{f}^a$ is the advective flux and $ \blocktensor{f}^{\nu}$ is the diffusive flux.
Because of the dependency of the diffusive flux on $\Nabla \mathbf{q}$, \eqref{eq:NScons} is a second order PDE.
Following Arnold et al. \cite{Arnold2002}, \eqref{eq:NScons} can be rewritten as a first-order system,
\begin{subequations}\label{eq:SystemAdvDiff}
\begin{empheq}[left=\empheqlbrace]{align} 
\partial_t \mathbf{q} 
+ \Nabla \cdot \left( \blocktensor{f}^a (\mathbf{q}) - \blocktensor{f}^{\nu} (\mathbf{q}, \blocktensor{g}) \right)&= \mathbf{0} \ , \ \text{in } \Omega,  \label{eq:OuterEq} \\
\Nabla \mathbf{q}
&=  \blocktensor{g}  , \  \text{in } \Omega. \label{eq:InnerEq}
\end{empheq}
\end{subequations}

To obtain the DGSEM-version of \eqref{eq:SystemAdvDiff}, the computational domain is subdivided into non-overlapping hexahedral elements, all variables are approximated by piece-wise Lagrange interpolating polynomials of degree $N$ that are continuous in each element, but allowed to be discontinuous across element interfaces: $\mathbf{q} \leftarrow \mathbf{q}^N$, $\blocktensor{f} \leftarrow \blocktensor{f}^N$ and $\blocktensor{g} \leftarrow \blocktensor{g}^N$.
Furthermore, \eqref{eq:OuterEq} and \eqref{eq:InnerEq} are multiplied by an arbitrary polynomial (test function) of degree $N$, the derivative terms are integrated by parts, and all integrals are evaluated numerically with a quadrature rule of $N+1$ points, to obtain
\begin{subequations}\label{eq:DGSEMsystem}
\begin{empheq}[left=\empheqlbrace]{align} 
J_j w_j \partial_t \mathbf{q}^N_j 
- \int_{\Omega^e}^N {\blocktensor{f}}^N \cdot \Nabla {\phi_j} \d \Omega^e 
+ \int_{\partial \Omega^e}^N \numflux{f} {\phi_j} \textrm{d} S^e 
&= \mathbf{0}, \label{eq:DGSEMsystem:1}\\
- \int _{\Omega^e}^N \mathbf{q}^N \Nabla \phi_j \d \Omega^e 
+ \int_{\partial \Omega^e}^N \phi_j \hat{\mathbf{q}} \vec{n} \d S^e
&=
J_j w_j  \blocktensor{g}^N_j \label{eq:DGSEMsystem:2}
\end{empheq}
\end{subequations}
for each degree of freedom of each element.
In \eqref{eq:DGSEMsystem}, $\numflux{f}$ and $\numflux{q}$ are the numerical traces of the flux and the solution, respectively, the functions $\phi_j$ are the so-called basis functions, which are tensor product expansions of the Lagrange interpolating polynomials, the $J_j$ are the Jacobians of the geometry transformation with which the mesh is created, and the $w_j$ are the weights of the quadrature rule. 
The derivation of \eqref{eq:DGSEMsystem} is given in \cite{kopriva2009implementing,Gassner2009a}.\\ 

The discretization of the system can be compactly written as
\begin{equation}\label{eq:DGMassMatrixLeft}
\mat{M}^N \bigderiv{\mathbf{Q}^N}{t} + \F^N(\mathbf{Q}^N) = \mathbf{0},
\end{equation}
where $\mat{M}^N$ is the mass matrix of the system, $\mathbf{Q}^N$ a vector with all the unknowns, and $\F^N$ the discrete partial differential operator. 
The mass matrix of the DGSEM is diagonal so, \eqref{eq:DGMassMatrixLeft} is often rewritten as
\begin{equation}\label{eq:DGMassMatrixRight}
\bigderiv{\mathbf{Q}^N}{t} + (\mat{M}^N)^{-1}\F^N(\mathbf{Q}^N) = \mathbf{0}.
\end{equation}

\subsection{Formulation of the Truncation Error} \label{sec:Uadapt:NewTruncError}

The truncation error is defined as the difference between the discrete partial differential operator and the continuous partial differential operator, both applied to the exact solution of the problem. 
This is often known as Generalized Truncation Error Expresion (GTEE) \cite{roy2009strategies,oberkampf2010verification}. 
For the problem at hand, we take the difference between the discrete equation \eqref{eq:DGMassMatrixRight} applied to the sampled continuous solution, $\mathbf{I}^N \mathbf{q}$, and the sampled continuous equation \eqref{eq:NScons} to obtain:
\begin{equation}\label{eq:TEderivation1}
\bigderiv{\mathbf{I}^N \mathbf{q}}{t} - \mathbf{I}^N \partial_t \mathbf{q} + (\mat{M}^N)^{-1}\F^N(\mathbf{I}^N \mathbf{q}) -  \mathbf{I}^N F(\mathbf{q})= \tilde{\tauvec}^N,
\end{equation}
where $ \tilde{\tauvec}^N$ is the truncation error and $\mathbf{I}^N$ is a restriction/prolongation operator used to project the solution from one space to another. Here we are simply sampling the continuous solution into our discrete space. 

Assuming that the restriction operator commutes with the time derivative ($\mathbf{I}^N \partial_t \mathbf{q}=\bigderiv{\mathbf{I}^N \mathbf{q}}{t}$), the two first terms in \eqref{eq:TEderivation1} cancel out, \begin{equation}\label{eq:TEtildedef}
\tilde{\tauvec}^N = (\mat{M}^N)^{-1}\F^N(\mathbf{I}^N \mathbf{q}) -  \mathbf{I}^N F(\mathbf{q}).\\
\end{equation}

The discretization error is defined as the difference between the approximate solution and the exact solution to the problem, i.e., $\epsvec^N=\mathbf{Q}^N-\mathbf{I}^N \mathbf{q}$. Taking the difference between \eqref{eq:DGMassMatrixRight} and \eqref{eq:NScons} sampled into the discrete space, we get the following.
\begin{equation}\label{eq:DETE1}
\bigderiv{\mathbf{Q}^N}{t} - \mathbf{I}^N \partial_t \mathbf{q} + (\mat{M}^N)^{-1}\F^N(\mathbf{Q}^N) -  \mathbf{I}^N F(\mathbf{q})= \mathbf{0}.
\end{equation}
By the linearity of the time derivative and using the definition of the discretization error,
\begin{equation}\label{eq:DETE2}
\bigderiv{\epsvec^N}{t} = - (\mat{M}^N)^{-1}\F^N(\epsvec^N+\mathbf{I}^N \mathbf{q}) +  \mathbf{I}^N F(\mathbf{q}).
\end{equation}
Now, for simplicity, we consider that $\F^N$ is a linear operator (it can be linearized otherwise to achieve a similar result; see, e.g., \cite{tyson2019relinearization}),
\begin{equation}\label{eq:DETE3}
\bigderiv{\epsvec^N}{t} + (\mat{M}^N)^{-1}\F^N(\epsvec^N)= - (\mat{M}^N)^{-1}\F^N(\mathbf{I}^N \mathbf{q}) +  \mathbf{I}^N F(\mathbf{q})=-\tilde{\tauvec}^N.
\end{equation}
Therefore, for linear operators, the discretization error is governed by the same equation as the numerical solution with the addition of the truncation error as a source term. This equation is known as Discrete Error Transport Equation (DETE). As can be seen, the truncation error is interesting for mesh adaptation, as it acts as a source for the generation of discretization error.\\

The truncation error definition, \eqref{eq:TEtildedef},  can be re-scaled with the mass matrix:
\begin{equation} \label{eq:TEdef}
\tauvec^N = \mat{M}^N \tilde{\tauvec}^N = \F^N(\mathbf{I}^N {\mathbf{q}}) - \mat{M}^N \mathbf{I}^N F({\mathbf{q}}),
\end{equation}
which is equivalent to defining the truncation error as the projection of the difference between discrete and continuous operators on the individual basis functions of the finite element subspace. 
With this definition, the DETE reads:
\begin{equation}\label{eq:DETE4}
\mat{M}^N \bigderiv{\epsvec^N}{t} + \F^N(\epsvec^N)= - \F^N(\mathbf{I}^N \mathbf{q}) +  \mat{M}^N \mathbf{I}^N F(\mathbf{q})=-\tauvec^N.
\end{equation}

The definition \eqref{eq:TEdef} has previously been used in \cite{Rubio2015,Kompenhans2016,Kompenhans2016a,RuedaRamirez2019,RuedaRamirez2019a}, and will be called traditional formulation in this work. The definition \eqref{eq:TEtildedef} has been recently used in \cite{laskowski2022functional} and will be called, in this work, new formulation. 
The main advantage of the new formulation is that it is directly related to functional errors, as shown in \cite{laskowski2022functional}.

Note that previous works in the context of the DGSEM \cite{Rubio2015,Kompenhans2016,Kompenhans2016a,RuedaRamirez2019,RuedaRamirez2019a,laskowski2022functional} focused primarily on steady-state problems; therefore, the second term in the RHS of \eqref{eq:TEtildedef} and \eqref{eq:TEdef} was zero.\\

We have derived two formulations of the truncation error and showed that each leads to a different version of the DETE: \eqref{eq:DETE3} and \eqref{eq:DETE4}. 
Although both formulations appear to be similar, some remarks can be made about their properties.

\begin{enumerate}
\item The traditional version of the truncation error \eqref{eq:TEdef} acts as a source term for the discretization error, after projecting it point-wise on the basis functions $\phi_j$, which build the finite element subspace, as shown in \eqref{eq:DETE4}.
On the other hand, the new approximation of the truncation error \eqref{eq:TEtildedef} acts directly as a source term of the pointwise values of the discretization error; see \eqref{eq:DETE3}. 

\item Since DGSEM is a collocation method and the mass matrix is a diagonal matrix containing the mapping Jacobian and quadrature weights, see, for example, \cite{ferrer2022horses3d}, each form of the truncation error can be obtained from the other by scaling it point-wise with $J_j w_j$, see \eqref{eq:TEdef}.
In other words, the main difference between both formulations is the weight that they give to the element size.

\item Due to the strong similarities between the two truncation error approximations, the anisotropic properties and the possibility to estimate the error in a multigrid cycle (see \cite{RuedaRamirez2019a}) hold for both error measures.

\end{enumerate}

Due to the similarities of both truncation error formulations, the tilde notation will be dropped in next sections, and the expressions will hold for both, unless the contrary is explicitly stated. 

On a final note, in \cite{Rubio2015} the concept of \emph{isolated} truncation error was introduced. 
While the standard \emph{non-isolated} truncation error uses all terms appearing in the discrete discontinuous Galerkin variational formulation, \eqref{eq:DGSEMsystem}, the \emph{isolated} truncation error replaces the surface numerical flux functions by simple evaluations of the flux with the inner solution of each element. 
The \emph{isolated} truncation error has some advantages for mesh adaptation, as shown in \cite{Kompenhans2016a}. 

\subsection{Truncation Error Estimation in Unsteady Problems} \label{sec:Uadapt:Estim}

Now that we have defined the truncation error, we need a method to estimate it when the exact solution is not available.
Previous works, see \cite{Rubio2015, Kompenhans2016,Kompenhans2016a,RuedaRamirez2019,RuedaRamirez2019a,laskowski2022functional}, consider only steady problems, and, therefore, only the first term of \eqref{eq:TEtildedef} or \eqref{eq:TEdef} takes nonzero values. 
In these works, the exact solution is approximated by a solution obtained on a higher-order mesh ($P > N$),
\begin{equation}
\mathbf{I}^N \mathbf{q} \approx \mathbf{I}^N \mathbf{q}^P,
\end{equation}
in a process known as $\tau$-estimation \cite{Rubio2013,Rubio2015}. Using this solution, the truncation error is estimated in all coarser meshes $N\in[1,P-1]$. This is useful for adaptation, as the exact polynomial degree required for a given accuracy can be directly read from the estimated truncation errors. If the problem solved has spatial dimension higher than one, the coarser meshes can be generated with anisotropic polynomial degrees (different polynomial degrees for the different spatial dimensions), resulting in the so-called truncation error map. Additionally, the truncation error for polynomial degrees $N>P-1$ can be estimated by $\tau$-extrapolation, and the whole estimation process can be embedded within an anisotropic multigrid cycle. The interested reader is referred to \cite{RuedaRamirez2019a} and references therein for details. \\

In unsteady problems, we can either find an approximation for the second term of \eqref{eq:TEtildedef} and \eqref{eq:TEdef} or reduce the problem to a steady-state case. The two main alternatives are listed below.

\begin{enumerate}

\item \textbf{Variational DG form}: 
To find an approximation for the second term of \eqref{eq:TEtildedef} or \eqref{eq:TEdef}, it is also reasonable to rely on the high-order solution, $\mathbf{q}^P$, and to approximate the continuous partial differential operator by the high-order one, 
\begin{equation}
\mathbf{I}^N F(\mathbf{q}) \approx 
    \mathbf{I}^N (\mat{M}^P)^{-1}\F^P(\mathbf{q}^P).
\end{equation}

\item \textbf{Dual time-stepping}: The original PDE, \eqref{eq:NScons}, can be reformulated using a dual time-stepping technique \cite{Arnone1995,Rumsey1995} as
\begin{equation}
\frac {\partial \mathbf{q}}{\partial \tau} + \frac {\partial \mathbf{q}}{\partial t} + F(\mathbf{q}) = \mathbf{0},
\end{equation}
where $\tau$ is a pseudo-time that is marched to steady-state in every time step of the physical time, $t$.
In dual time-stepping methods, the time derivative is usually discretized with an implicit method, and the pseudo-time derivative is either discretized with implicit or explicit methods.

The advantage of using an explicit scheme for the pseudo-time derivative is that the physical problem is integrated in time implicitly, without having to solve linear systems. 
All in all, the discretized/sampled system yields
\begin{equation}
\mat{M}^N \frac{\delta_s}{\delta \tau} \mathbf{Q}^N + \mat{M}^N \frac{\delta_p}{\delta t} \mathbf{Q}^N + \F^N(\mathbf{Q}^N) = \mathbf{0},
\end{equation}
where $\delta_s / \delta \tau$ is an operator that imposes the chosen pseudo-time-integration scheme and $\delta_p / \delta t$ is the operator for the chosen time-integration scheme.

Since we now have a steady problem, it is possible to consider only the first term of \eqref{eq:TEtildedef} or \eqref{eq:TEdef} if the partial differential operators are redefined as
\begin{equation}
 \Res^N (\mathbf{Q}^N)
= \mat{M}^N \frac{\delta_p}{\delta t} \mathbf{Q}^N + \F^N(\mathbf{Q}^N),
\end{equation}
for the traditional formulation, and
\begin{equation}
\tilde{\Res}^N (\mathbf{Q}^N)
= \frac{\delta_p}{\delta t} \mathbf{Q}^N + \mat{M}^{-1} \F^N(\mathbf{Q}^N),
\end{equation}
for the new formulation.

\end{enumerate}

The dual time-stepping truncation error estimation approach is useful to implement in codes that already use a dual time-stepping integration method.
The variational DG approach is, in general, simpler to implement and imposes almost no overhead, since the term $(\mat{M}^P)^{-1}\F^P(\mathbf{q}^P)$ is obtained in the calculation of $\mathbf{q}^P$, required for the estimation of the truncation error.
In the rest of this paper, we use the variational DG form of the unsteady truncation error.

\subsection{$p$-Adaptation Strategies} \label{sec:Uadapt:pAdaptAlg}

Two adaptation strategies can be identified in unsteady flow simulations: dynamic and static adaptation. 
These two strategies have already been widely used for unsteady adaptivity. 
See, for example, \cite{Blaise2012,Cagnone2012} for dynamic adaptation methods or \cite{Naddei2018,Fidkowski2011} for static adaptation methods.

In the following sections, we present a detailed description of how dynamic and static $p$-adaptation can be implemented for truncation error-based $p$-adaptation methods.
The main difference between the strategies presented here and those encountered in the literature is the way the adaptation algorithms treat the error estimates.

Most of the error estimation strategies available in the literature are designed to mark a number of elements for enrichment or order reduction \cite{Naddei2018,Cagnone2012,Blaise2012}.
As a result, in every adaptation stage, the polynomial degrees are increased or reduced by one.
On the contrary, the truncation error estimation provides an exact value for the polynomial degree needed for each coordinate direction of every element after each estimation stage. 
This property of the $\tau$-estimation method provides several advantages for the $p$-adaptation of unsteady computations, as will be discussed in the following sections.

\subsubsection{Dynamic $p$-Adaptation} \label{sec:DynamicPAdapt}

The dynamic $p$-adaptation is the most straightforward $p$-adaptation strategy for unsteady flows.
It computes an error measure periodically during a simulation and adapts the polynomial degrees of the discretization according to the estimated error, right after every estimation procedure. The adaptation of the polynomial degree in every step follows the procedure introduced in \cite{RuedaRamirez2019a}.
The dynamic $p$-adaptation strategy is well suited for transient simulations in which the region of interest changes over time. 

\begin{figure}
\centering
\begin{tikzpicture}
\node (start) [startbox] {Initialize};
\node (advance) [process, below of=start, yshift=-0.5cm] {Advance $\Delta t_e$ in time};
\draw [arrow] (start) -- node[name=y]{} (advance);
\node (dec1) [decision, below of=advance, yshift=-0.5cm] {Reached final time?};

\node (estimate) [process, right of=dec1, xshift = 4cm] {Estimate error};
\node (adapt) [process, right of=y, xshift = 4cm] {\uline{Adapt polynomial degrees}};
\node (end) [endbox, below of=dec1,yshift = -0.5cm]{Finalize};
\draw [arrow] (advance) -- (dec1);

\draw [arrow] (dec1) -- node[anchor=south] {no} (estimate);
\draw [arrow] (estimate) -- (adapt);
\draw [arrow] (adapt) -- (y);
\draw [arrow] (dec1) -- node[anchor=east] {yes} (end);
\end{tikzpicture}
\caption{Flowchart of the dynamic $p$-adaptation.}
\label{fig:FlowCDynamicPAdapt}
\end{figure}
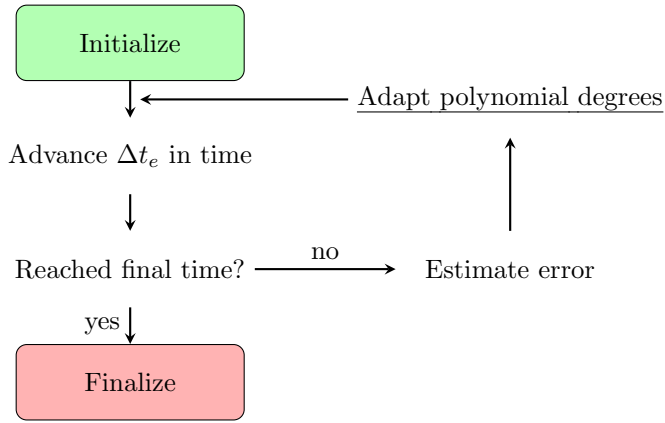

Figure \ref{fig:FlowCDynamicPAdapt} illustrates the dynamic $p$-adaptation process. 
The interval between adaptation stages, $\Delta t_e$, can be specified as a physical time, as a number of iterations (time steps) or can be changed throughout the simulation.
At every $p$-adaptation stage, the underlined process in Figure \ref{fig:FlowCDynamicPAdapt}, the storage must be reallocated and the solution projected in the new polynomial spaces.
Since this process is done several times during the solution procedure, some overhead is expected.
Therefore, the construction of the data structures for the new spatial resolution and the transfer of information are critical steps that must be optimized to enhance the performance.
Furthermore, if the $\tau$-estimation method is used, a number of low-order ($N<P$) discretizations are needed to evaluate the truncation error.
As a result, an additional overhead is added in the construction of these coarse grids.

Note that a traditional error estimator, which simply marks some elements for refinement or coarsening, may perform poorly with the dynamic $p$-adaptation strategy of Figure \ref{fig:FlowCDynamicPAdapt}.
Such an error estimator imposes a one-by-one increase in the polynomial degree.
Therefore, if $\Delta t_e$ is too large, a dynamic $p$-adaptation strategy may not have enough time to increase the resolution of a zone of the domain before the flow feature of interest goes out of it.
In other words, the refinement zones are likely to lag behind the difficult-to-capture flow features.
On the contrary, since the truncation error estimator identifies what polynomial degree is needed immediately, the resolution can be increased to the necessary level immediately.
As a result, the truncation error estimator may be more suitable to handle larger values of $\Delta t_e$ than traditional estimators.

It is possible to obtain overshoots in the truncation error estimates if the polynomial degree of the reference mesh, $P$, is too low, as shown in \cite{RuedaRamirez2019a}.
Therefore, a truncation error-based dynamic $p$-adaptation method may suffer unnecessary polynomial degree oscillations that are caused and nurtured by the constant jump between a low and a high $P$. 
These polynomial degree oscillations may deteriorate the accuracy and, therefore, should be avoided when possible.
A possible way to attenuate this phenomenon is to limit the maximum polynomial degree jump (by element and coordinate direction) after each $p$-adaptation stage.

Additionally, in parallelized simulations, a dynamic $p$-adaptation strategy requires dynamic load balancing to maintain an even workload between the processors and avoid deadlocks.
Otherwise, the reduction in the number of degrees of freedom achieved with the enhanced spatial discretization may not translate into shorter computation times.
The design of efficient dynamic load balancing algorithms is a challenging topic of research that is not treated in this work.

\subsubsection{Static $p$-Adaptation}

The static $p$-adaptation strategy differs from the dynamic $p$-adaptation in that only one $p$-adaptation process is performed after several estimation stages.

\begin{figure}
\centering
\begin{tikzpicture}
\node (start) [startbox] {Initialize};
\node (advance) [process, below of=start, yshift=-0.5cm] {Advance $\Delta t_e$ in time};

\draw [arrow] (start) -- node[name=y]{} (advance);

\node (dec1) [decision, below of=advance, yshift=-0.5cm] {Reached final estimation time?};

\node (estimate) [process, right of=advance, xshift = 4cm] {Estimate error and save};
\node (adapt) [process, below of=dec1, yshift = -0.5cm] {\uline{Adapt polynomial degrees}};

\node (advance2) [process, below of=adapt, yshift = -0.5cm] {Advance in time until final time};

\node (end) [endbox, below of=advance2,yshift = -0.5cm]{Finalize};
\draw [arrow] (advance) -- (dec1);

\draw [arrow] (dec1) -| node[anchor=north] {no} (estimate);
\draw [arrow] (estimate) |- (y);
\draw [arrow] (dec1) -- node[anchor=east] {yes} (adapt);
\draw [arrow] (adapt) -- (advance2);
\draw [arrow] (advance2) -- (end);
\end{tikzpicture}
\caption{Flowchart of the static $p$-adaptation.}
\label{fig:FlowCStaticPAdapt}
\end{figure}
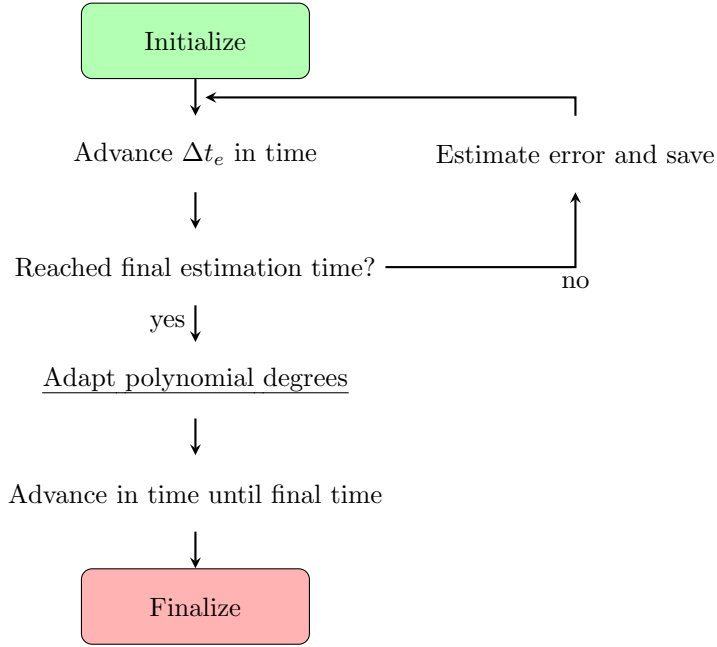

Figure \ref{fig:FlowCStaticPAdapt} presents a flowchart of the static adaptation strategy.
First, the solution is advanced in time with a fixed spatial resolution until a final estimation time, $T_e$, is reached. 
During this stage, periodic error estimations are performed with an interval of $\Delta t_e$ but, instead of changing the spatial resolution, the error estimation is stored for future processing.
When $T_e$ is reached, a $p$-adaptation procedure is performed using all the stored error estimates.
Subsequently, the simulation is advanced in time with the new fixed spatial resolution until the final time.
The static $p$-adaptation is well suited for simulations in which the features of interest are located in a fixed region of the domain. 

Note that a traditional error estimator, which only marks some elements for one-by-one refinement or coarsening, may also perform poorly in the static $p$-adaptation algorithm of Figure \ref{fig:FlowCStaticPAdapt}.
In fact, if such an error estimator is used, the algorithm would have to be slightly modified, so that after the $p$-adaptation stage, the simulation goes back to the error estimation stage, as in \cite{Naddei2018}. 
That extra loop would have to be repeated a specific number of times, or until no element is marked for refinement or coarsening, which represents increased computational cost. 
Therefore, the ability to predict the exact polynomial that is needed makes the truncation error estimator an attractive indicator for statically $p$-adapted unsteady simulations.

The static $p$-adaptation strategy provides several implementation advantages over the dynamic $p$-adaptation.
First, the construction of the data structures and the projection of the solution to the new spatial discretization are no longer critical steps, as the $p$-adaptation procedure is only done once.
Therefore, these operations can even be performed off-line, and their computational cost does not significantly impact the performance of the method.
Second, the coarse-grid discretizations that are needed for the truncation error estimation are only constructed once at the beginning of the simulation and used throughout the entire estimation stage.
Finally, dynamic load balancing is no longer needed, as the loads must be balanced only once after the $p$-adaptation step.

The static $p$-adaptation algorithm has two drawbacks.
First, it is only useful for statistically steady flows or where the features that need high spatial resolution are located in a specific region of the domain.
Second, the static $p$-adaptation strategy requires a preliminary simulation to estimate the error. 

In the case of aerodynamic simulations of external flow, the flows are usually statistically steady, and the \textit{interesting} flow features are concentrated in a small region of the domain.
Furthermore, the estimation time is generally much shorter than the total simulation time.\\

To process the $p$-anisotropic truncation error estimation data and feed the $p$-adaptation algorithm, two main approaches can be identified:

\begin{enumerate}
\item At each estimation stage, $s$, select the polynomial degrees for each of the elements in the mesh, $N_i^{e,s}$, and then predict a final polynomial degree from the estimates, 
\begin{equation}
N_i^e = F(N_i^{e,1}, \dots, N_i^{e,n_e}),
\end{equation}
where $n_e$ is the number of estimation stages. 

\item At each estimation stage, $s$, generate the truncation error map (see \cite{RuedaRamirez2019a} for details) for every element of the mesh, $\norm{\tauvec}_{\infty}^{e,s}$, compute a total truncation error map for every element based on the estimations, 
\begin{equation} \label{eq:Uadapt:totalTruncMap}
\norm{\tauvec}_{\infty}^e = F \left( \norm{\tauvec}_{\infty}^{e,1}, \dots, \norm{\tauvec}_{\infty}^{e,n_e} \right),
\end{equation}
and select the polynomial degree from the total truncation error map.
\end{enumerate}

The function $F(\cdot)$ can be defined in multiple ways, for example, the average or maximum functions, 
\begin{equation}
F_{\mathrm{av}} (N_i^{e,1}, \dots, N_i^{e,s}) = \frac{1}{n_e} \sum_{s=1}^{n_e} N_i^{e,s}, \ \ \ \ \ \
F_{\max} (N_i^{e,1}, \dots, N_i^{e,s}) =  \max_{s} |N_i^{e,s} |.
\end{equation}

As a conservative criterion, we use $F_{\max}$ to ensure that the specified truncation error threshold is satisfied throughout the whole simulation.\\

\textit{Approach 1} needs less storage space and can be implemented more easily than \textit{approach 2}. 
However, it may lead to the over-enrichment of some areas of the domain when combined with $p$-anisotropic discretizations. 
To illustrate this, let us consider a specific element in a hypothetical two-stage estimation procedure ($n_e = 2$) of a 2D simulation that uses $F_{\max}$. 
Let us assume a minimum polynomial degree $N_{\min}=1$  and a maximum polynomial degree $N_{\max}=3$ for the $p$-adaptation.
Furthermore, we are interested in selecting the polynomial degree combination that minimizes the number of degrees of freedom ($\NDOF$). 

Table \ref{tab:Uadapt:approach1} shows a possible outcome of the two-stage estimation procedure.
There are five polynomial degrees that fulfill the specified error threshold, $\tau_{\max}$, in each estimation stage.
Among those, \textit{approach 1} would select the polynomial degrees underlined in red in each estimation stage because they minimize the instantaneous $\NDOF$.
As can be observed, at the end of the estimation, \textit{approach 1} selects the polynomial degrees $N^e=(3,3)$, which correspond to $\NDOF = 16$. 
This outcome is not optimal because $N^e=(2,2)$, with an associated $\NDOF = 8$, would actually fulfill $\tau_{\max}$ with fewer degrees of freedom.

\begin{table} [h]
\centering
\caption{Possible outcome of the \textit{approach 1} to static $p$-adaptation.} \label{tab:Uadapt:approach1}
\begin{tabular}{c|ccccc|ccccc|c}
Coordinate & \multicolumn{10}{c|}{Polynomial degrees with $\norm{\tauvec^N}_{\infty} < \tau_{\max}$} & Selected \\ 
direction & \multicolumn{5}{c|}{Stage $s=1$} & \multicolumn{5}{c|}{Stage $s=2$} & degree \\ 
$i$     & \multicolumn{5}{c|}{$N_i^{e,1}$} & \multicolumn{5}{c|}{$N_i^{e,2}$} & $N_i^e$ \\ \hline
$1$     & $3$ & $2$ & $3$ & $2$ & {\color{red} $\uline{1}$}  & $3$ & $2$ & $3$ & $2$ & {\color{red} $\uline{3}$}  & $3$ \\ 
$2$     & $3$ & $3$ & $2$ & $2$ & {\color{red} $\uline{3}$}  & $3$ & $3$ & $2$ & $2$ & {\color{red} $\uline{1}$}  & $3$ \\ \hline \hline
$\NDOF$ & $16$ & $12$ & $12$ & $9$ & {\color{red} $\uline{8}$}  & $16$ & $12$ & $12$ & $9$ & {\color{red} $\uline{8}$}  & $16$ 
\end{tabular}

\end{table}

\textit{Approach 2} generates a total truncation error map by applying $F_{\max}$. 
The use of $F_{\max}$ implies that the polynomial degree combinations that fulfill $\tau_{\max}$ in the total map are those that fulfill $\tau_{\max}$ in all the estimation stages, i.e. the combinations that are not underlined in red.
Of that set of combinations, \textit{approach 2} clearly selects $N^e=(2,2)$, with an associated $\NDOF=8$.

As shown in this simple example, \textit{approach 2} is better than \textit{approach 1}.
Therefore, \textit{approach 2} is selected for the static $p$-adaptation simulations that are shown in this paper.

Although superior, \textit{approach 2} has two drawbacks. 
First, the extrapolated truncation error map (see, e.g., \cite{RuedaRamirez2019a}) must be obtained for each estimation process, so that the total truncation error map can be obtained with \eqref{eq:Uadapt:totalTruncMap}, which involves more computational resources per estimation stage.
Second, if an element is not in the asymptotic range in any of its reference coordinate directions, it may not be possible to extrapolate the values of the inner truncation error map.
In such cases, instead of extrapolating the truncation error, we assign a high value to it for $N_i \ge P_i$ as a secure criterion.

\section{Numerical Results} \label{sec:Uadapt:Results}

In this section, we test the performance of the methods described in this paper to perform $p$-adaptation of unsteady flow problems using truncation error estimates.
The methodology presented in this paper is valid for the \textit{non-isolated} and \textit{isolated} truncation errors. For simplicity, we only use the \textit{isolated} truncation error in this section to drive the $p$-adaptation procedures. 

The $p$-adaptation procedures are implemented in the open source high-order discontinuous Galerkin framework HORSES3D \cite{ferrer2022horses3d}.
All simulations use the Roe solver \cite{Roe1981} as the advective numerical flux and BR1 \cite{Bassi1997} as the diffusive numerical flux.
The time-marching scheme in the following examples is Williamson's low-storage third-order Runge-Kutta method \cite{williamson1980low}.
Additionally, the time-step size is dynamically changed using the CFL condition in all simulations (see \cite{ferrer2022horses3d} for details).
The main reason is that we want to take as large time steps as possible, and the time-step size is a function, among others, of the polynomial degree.

Since the flow features that we analyze have a periodicity in time, the interval between estimation/adaptation stages, $\Delta t_e$, is selected as a constant time for each simulation and not as the time that corresponds to a number of time steps.

\subsection{Advection of a Density Pulse in a Uniform Flow} \label{sec:FlowVortex}

\begin{figure}[tb]
\begin{center}

\includegraphics[width=0.45\textwidth]{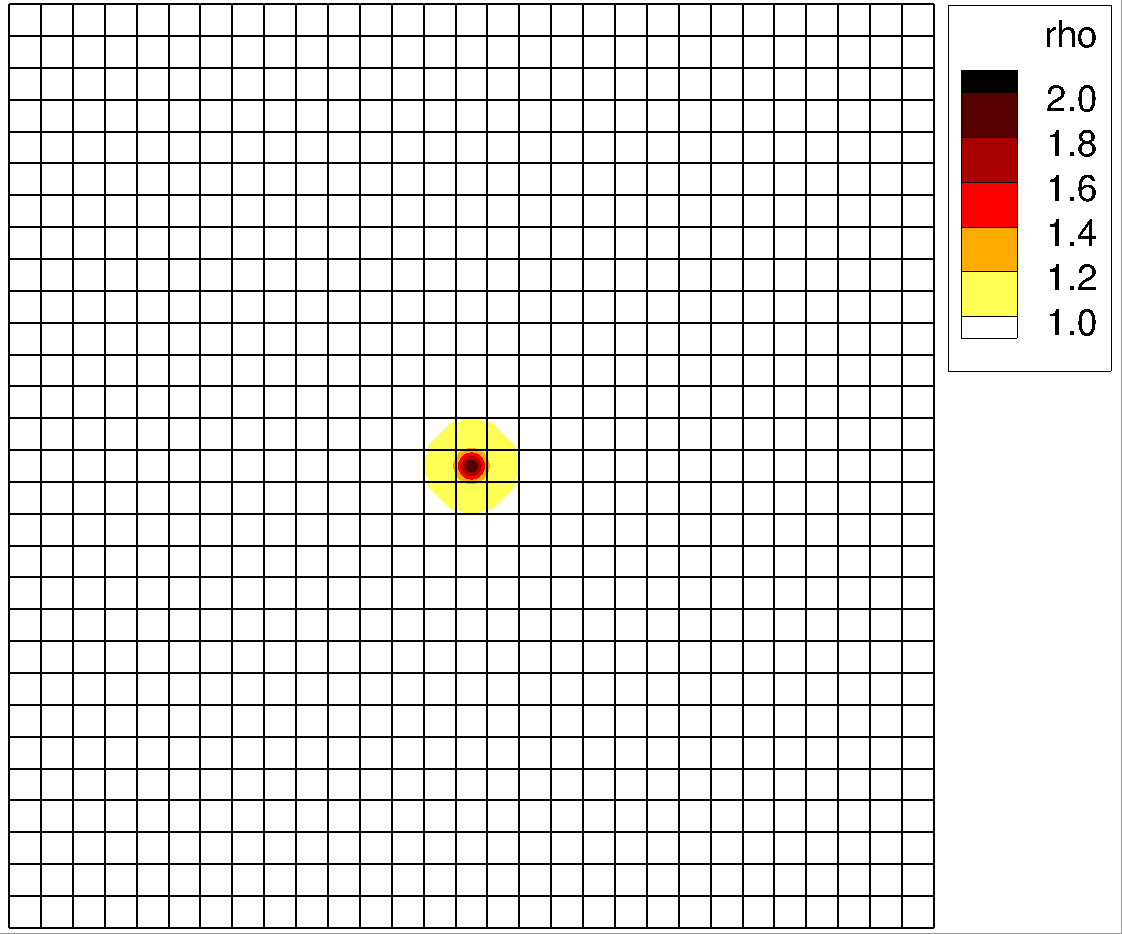} 
\caption{Initial condition of the advected density pulse simulation.}\label{fig:VortexContour}
\end{center}
\end{figure}

In this section, we simulate the advection of a Gaussian pulse in a square domain with periodic boundary conditions with the compressible Euler equations of gas dynamics and $\Ma=0.5$. 
The initial condition is
\begin{equation}
\begin{array}{ccl}
u = 1 & & \rho  = e^{5(x^2 + y^2)} + 1 \\
v = 0 & & p = 1
\end{array},
\end{equation}
and the two-dimensional computational domain is tessellated with a structured mesh of 841 quadrilateral elements, as shown in Figure \ref{fig:VortexContour}.
The final time is $t=29$, when the pulse should be back where it started. All simulations are run in serial with a sixth generation 8-core Intel i7 processor and 32GB of RAM.

In this test case, both forms of the truncation error (traditional and new) perform equivalently.
The polynomial degree distributions obtained with the new formulation of the truncation error are almost identical to those obtained with the traditional formulation, if the specified error threshold, $\tau_{\max}$, is scaled with the (constant) element size.
As we discussed in Section \ref{sec:Uadapt:NewTruncError}, the main difference between the two formulations is that the traditional truncation error is scaled with the element size. Since the element size is uniform in the whole domain, no significant differences are observed between the traditional and new forms of the truncation error.

We tested static and dynamic $p$-adaptation algorithms with truncation error thresholds ranging between $10^{-3} \le \tau_{\max} \le 1$, and intervals between adaptation/estimation stages ranging between $0.5 \le \Delta t_e \le 10$.
The polynomial degree is adapted according to the output of the error estimations in the range $1 \le N_i \le 8$ for each direction $i$ of every element.
Furthermore, the two different polynomial degree jump conditions that were introduced in \cite{RuedaRamirez2019a} are considered:

\begin{enumerate}[label=(\alph*)]
\item The first one imposes that the polynomial degree after every adaptation stage must fulfill
\begin{equation} \label{eq:Uadapt:N23}
N^{e}_i \ge \max_{j \in \mathcal{N}(e)} \bigg \lfloor \frac{2}{3} N^{j}_i \bigg \rfloor,
\end{equation}
where $N^{e}_i$ is the polynomial degree of element $e$ in the coordinate direction $i$, $\mathcal{N}(e)$ is the list of the neighbor elements of $e$, $N^{j}_i$ is the polynomial degree of the neighbor element $j$ in the matching coordinate direction $i$, and $\lfloor \cdot  \rfloor$ is the integer part floor function. 

\item The second polynomial degree jump condition imposes
\begin{equation} \label{eq:Uadapt:N_1}
N^{e}_i \ge \max_{j \in \mathcal{N}(e)} \bigg \lfloor N^{j}_i - 1 \bigg \rfloor.
\end{equation}
\end{enumerate}

In the static $p$-adaptation cases, a preliminary simulation must be run to estimate the error, as seen in Figure \ref{fig:FlowCStaticPAdapt}.
Since the pulse always changes position, the preliminary simulation must be run for $29$ time units to obtain a significant sample.
To have enough points to extrapolate the anisotropic truncation error estimates, the $\tau$-estimation simulation uses a discretization of uniform polynomial degree $P=3$.

\begin{figure}[h]
\begin{center}
\subfigure[Maximum error vs. number of DOFs.]{\label{fig:Vortex_N23_DOFs} \includegraphics[width=0.46\textwidth]{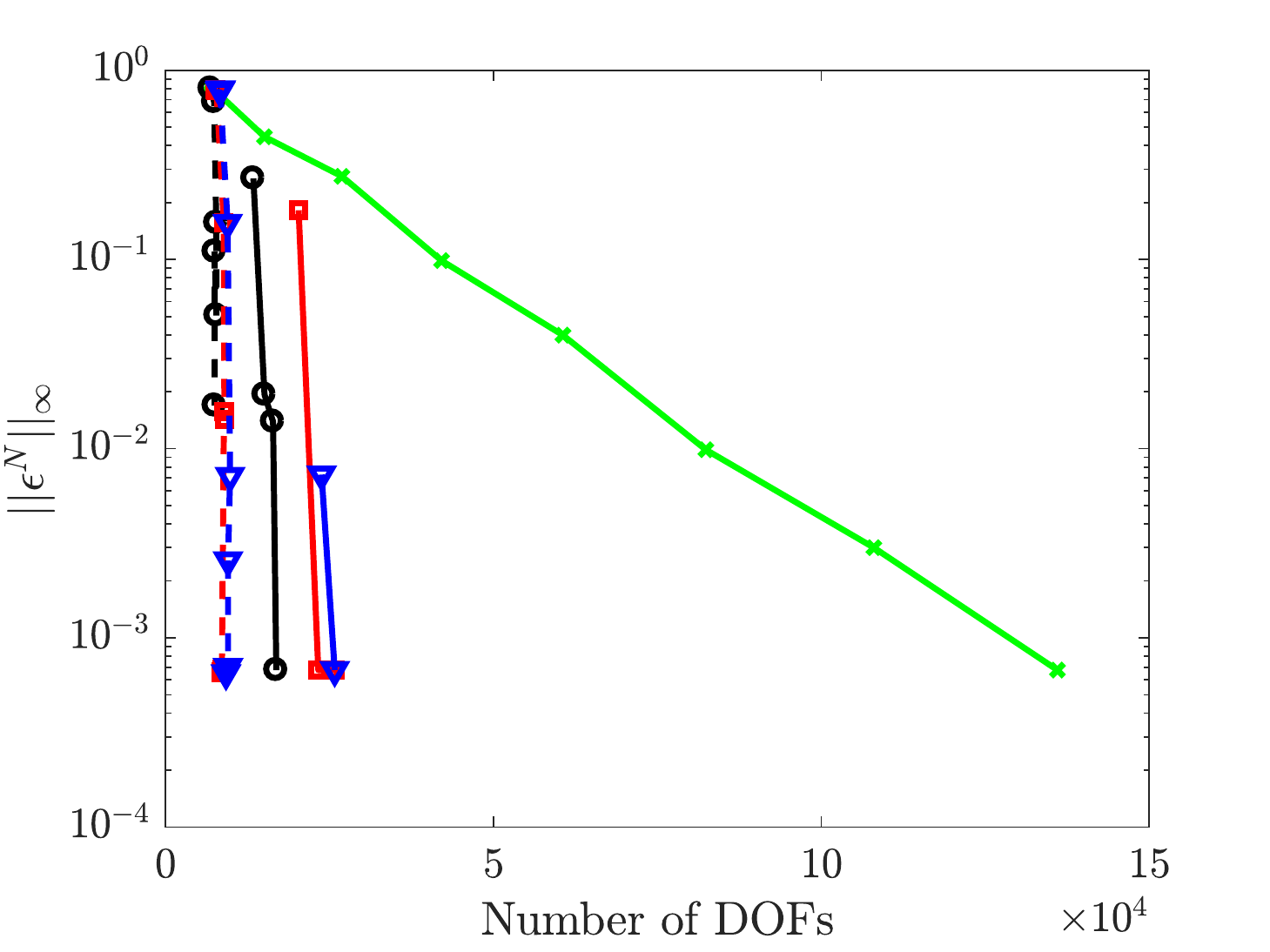}}\qquad
\subfigure[Maximum error vs. CPU-Time.]{\label{fig:Vortex_N23_Time} \includegraphics[width=0.46\textwidth]{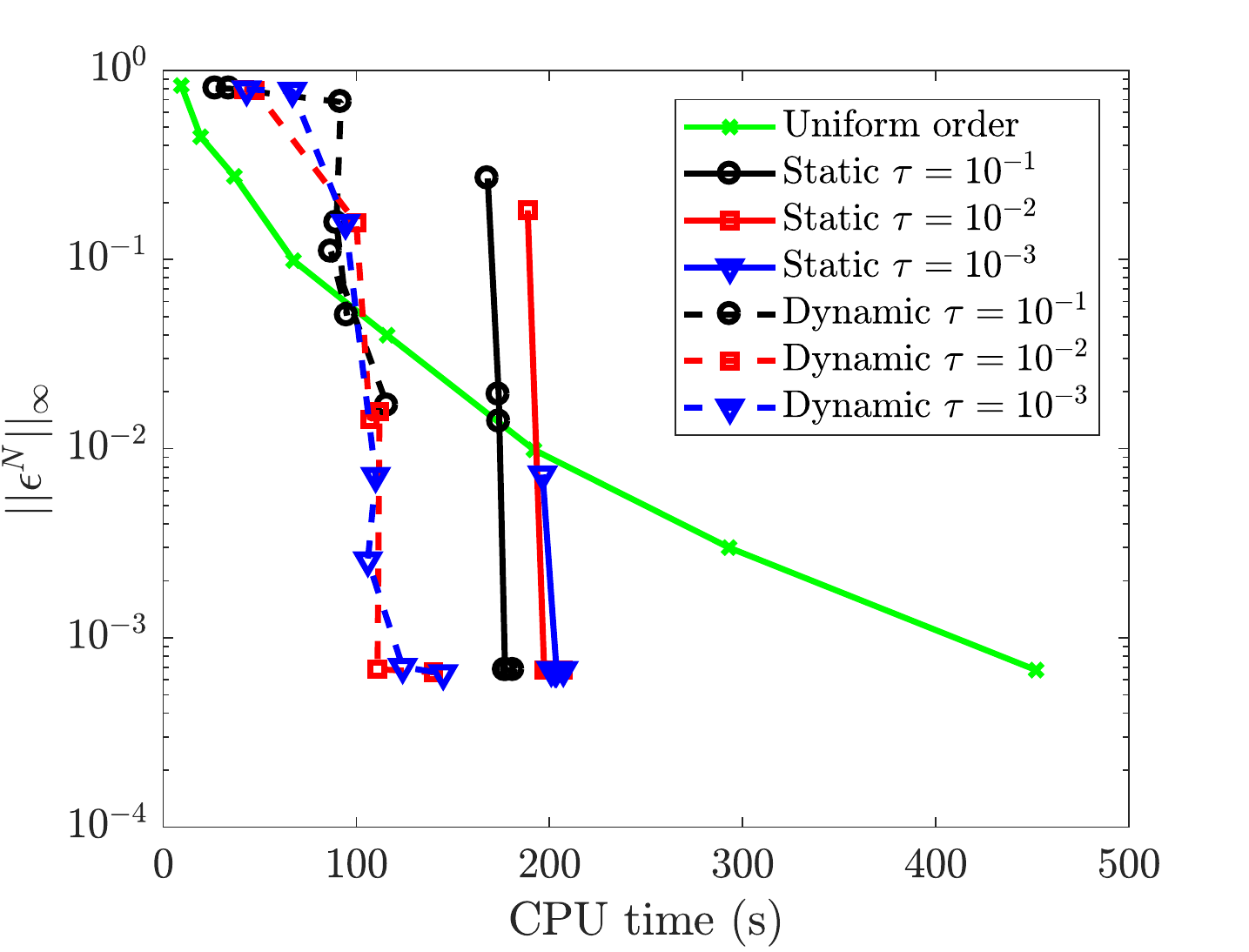}}

\caption{Error performance of static and dynamic $p$-adaptation procedures for the advection of a density pulse with a polynomial degree jump condition of $N^{e}_i \ge \max_{j \in \mathcal{N}(j)} \lfloor 2 N^{j}_i / 3 \rfloor$ \eqref{eq:Uadapt:N23}. 
Each point in the plot corresponds to an interval between adaptation/estimation stages, ranging between $0.5 \le \Delta t_e \le 10$. 
Different truncation error thresholds are represented with lines of different colors (black, red and blue).
Error performance of the uniform polynomial degree refinement (in green) is also included.  
 } \label{fig:VortexN23}
%
\subfigure[Maximum error vs. number of DOFs.]{\label{fig:Vortex_N_1_DOFs} \includegraphics[width=0.46\textwidth]{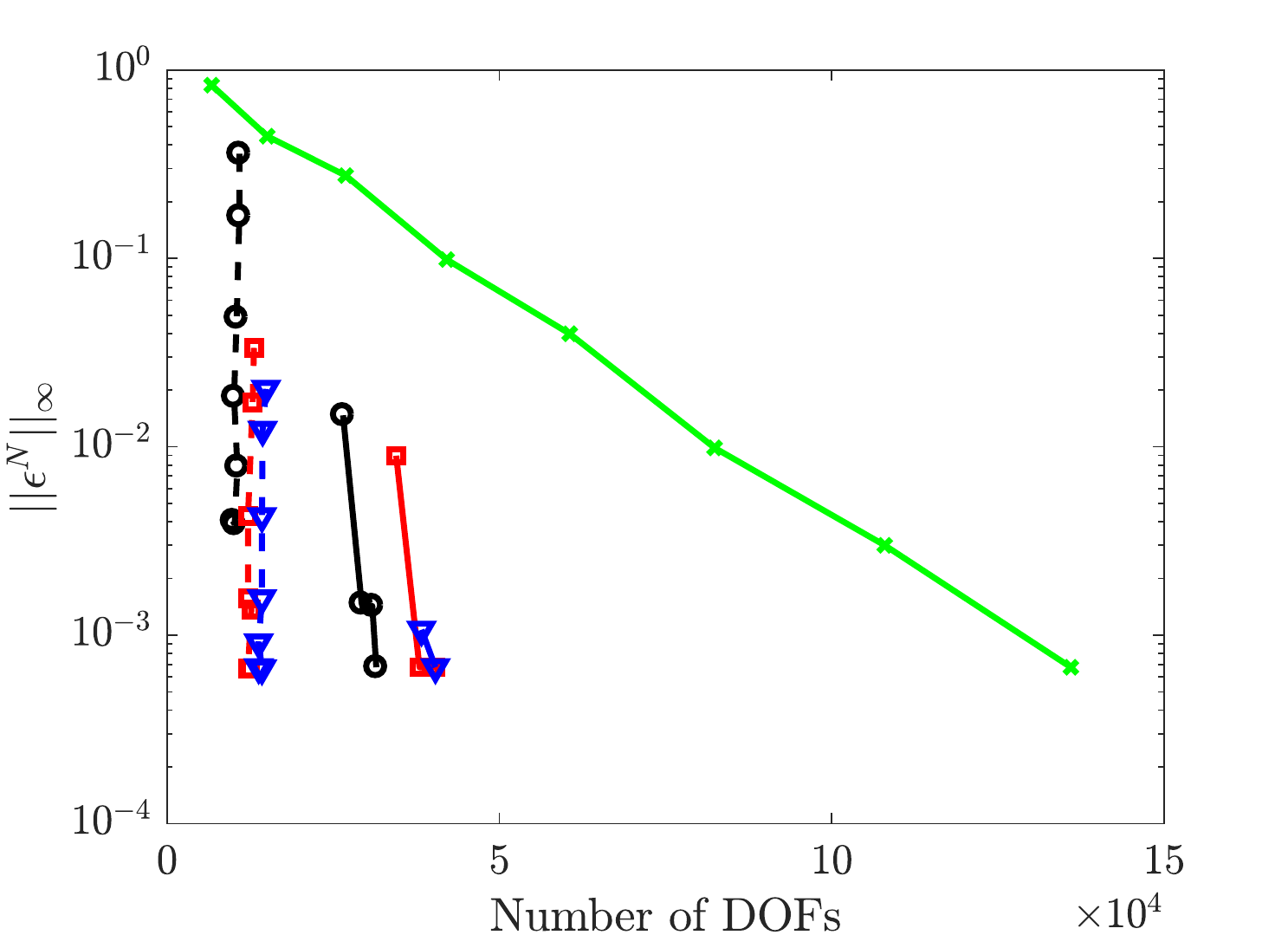}}\qquad
\subfigure[Maximum error vs. CPU-Time.]{\label{fig:Vortex_N_1_Time} \includegraphics[width=0.46\textwidth]{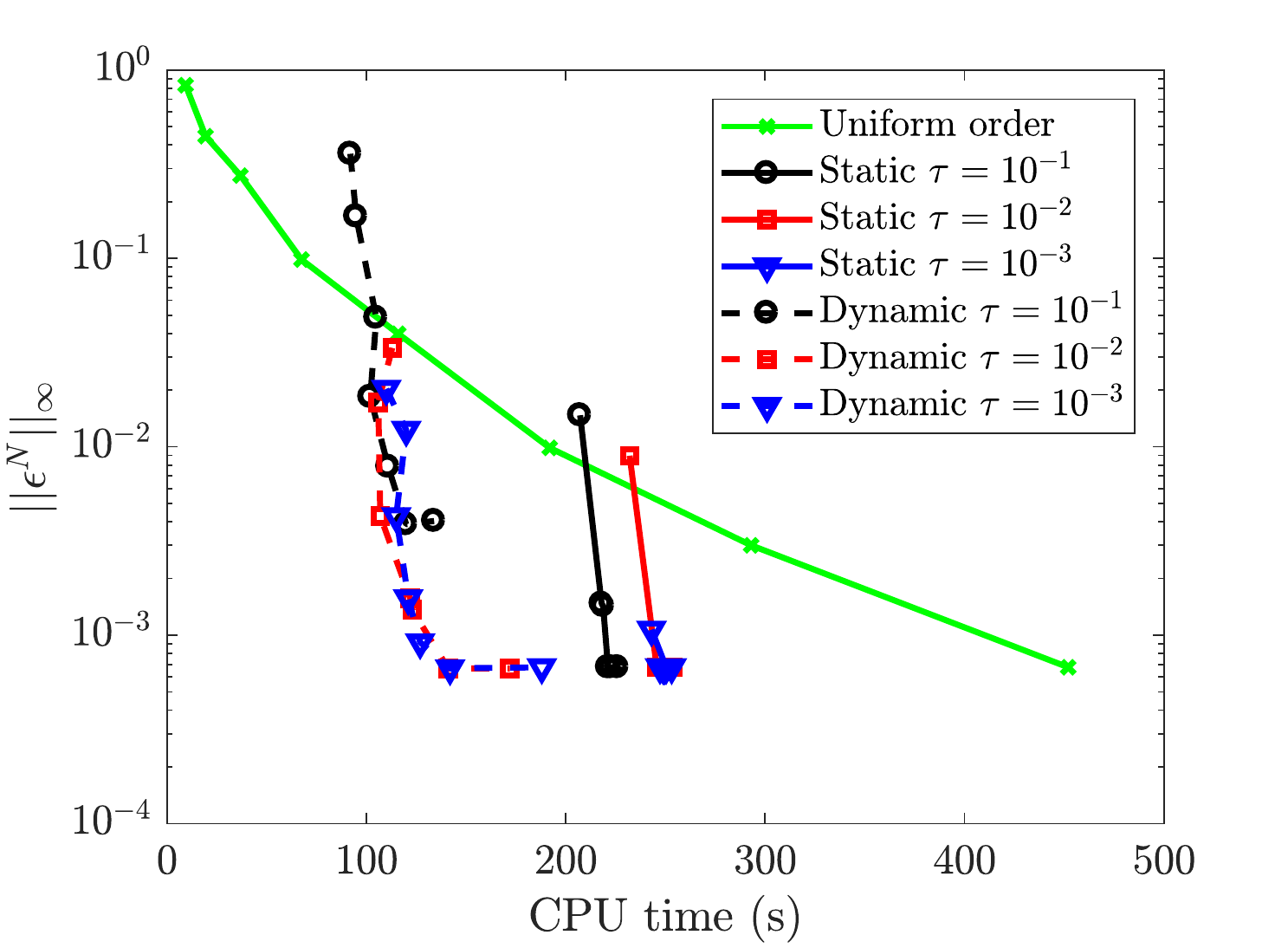}}

\caption{Error performance of static and dynamic $p$-adaptation procedures for the advection of a density pulse with a polynomial degree jump condition of $N^{e}_i \ge \max_{j \in \mathcal{N}(j)} \lfloor N^{j}_i - 1 \rfloor$ \eqref{eq:Uadapt:N_1}.
Each point in the graph corresponds to an interval between adaptation/estimation stages, ranging between $0.5 \le \Delta t_e \le 10$. 
Different truncation error thresholds are represented with lines of different colors (black, red and blue).
Error performance of the uniform polynomial degree refinement (in green) is also included. 
} \label{fig:VortexN_1}
\end{center}
\end{figure}

Figures \ref{fig:VortexN23} and \ref{fig:VortexN_1} show the behavior of the dissipation error as a function of the number of degrees of freedom and the computation time for the $p$-adaptive simulations with the polynomial degree jump conditions of \eqref{eq:Uadapt:N23} and \eqref{eq:Uadapt:N_1}, respectively.
The dissipation error is measured as the difference in $\rho$ between the exact solution and the simulation outcome at the centroid of the moving Gaussian.
The dispersion error, which can be measured as the absolute value of the position of the Gaussian centroid, is of the order of machine zero (the DGSEM exhibits very low dispersion errors in this case).
Note that the computation time needed for the estimation simulation in the static $p$-adaptation cases has already been added to the simulation time in Figures \ref{fig:Vortex_N23_Time} and \ref{fig:Vortex_N_1_Time}.

As can be observed, the truncation error-based $p$-adaptation techniques perform better than uniform refinement when a dissipation error $\norm{\gstate{\epsilon}^N}_{\infty} < 10^{-2}$ is desired, as they achieve the same errors with fewer degrees of freedom, which results in shorter computation times for a given accuracy.

Figures \ref{fig:Vortex_N23_DOFs} and \ref{fig:Vortex_N_1_DOFs} show that the number of degrees of freedom for the statically $p$-adaptive simulations is $2-4$ times higher than for the dynamically $p$-adaptive simulation.
This makes sense since the static $p$-adaptation algorithm enriches all the regions through which the pulse passes, whereas the dynamic $p$-adaptation algorithm effectively follows it.
The longer computation times that are observed in Figures \ref{fig:Vortex_N23_Time} and \ref{fig:Vortex_N_1_Time} for the statically $p$-adaptive simulations are not only the result of this effect, but also of the extra computation time invested in the preliminary $\tau$ estimation simulation.

\begin{figure}[h]
\begin{tabular}{cc}

\multicolumn{2}{c}{
	\includegraphics[width=0.45\textwidth]{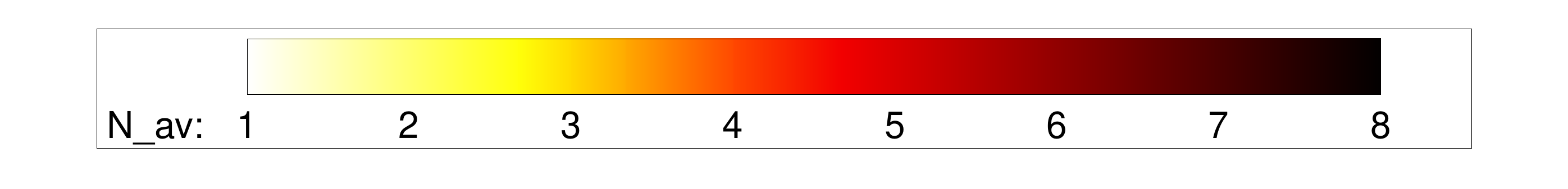}
} \\
\subfigure[$N^{+}_i \ge \lfloor N^{-}_i - 1   \rfloor$, $\Delta t_e = 10$.]
{
  \label{fig:VortexPol_N_1_10} 
  \includegraphics[width=0.45\textwidth]{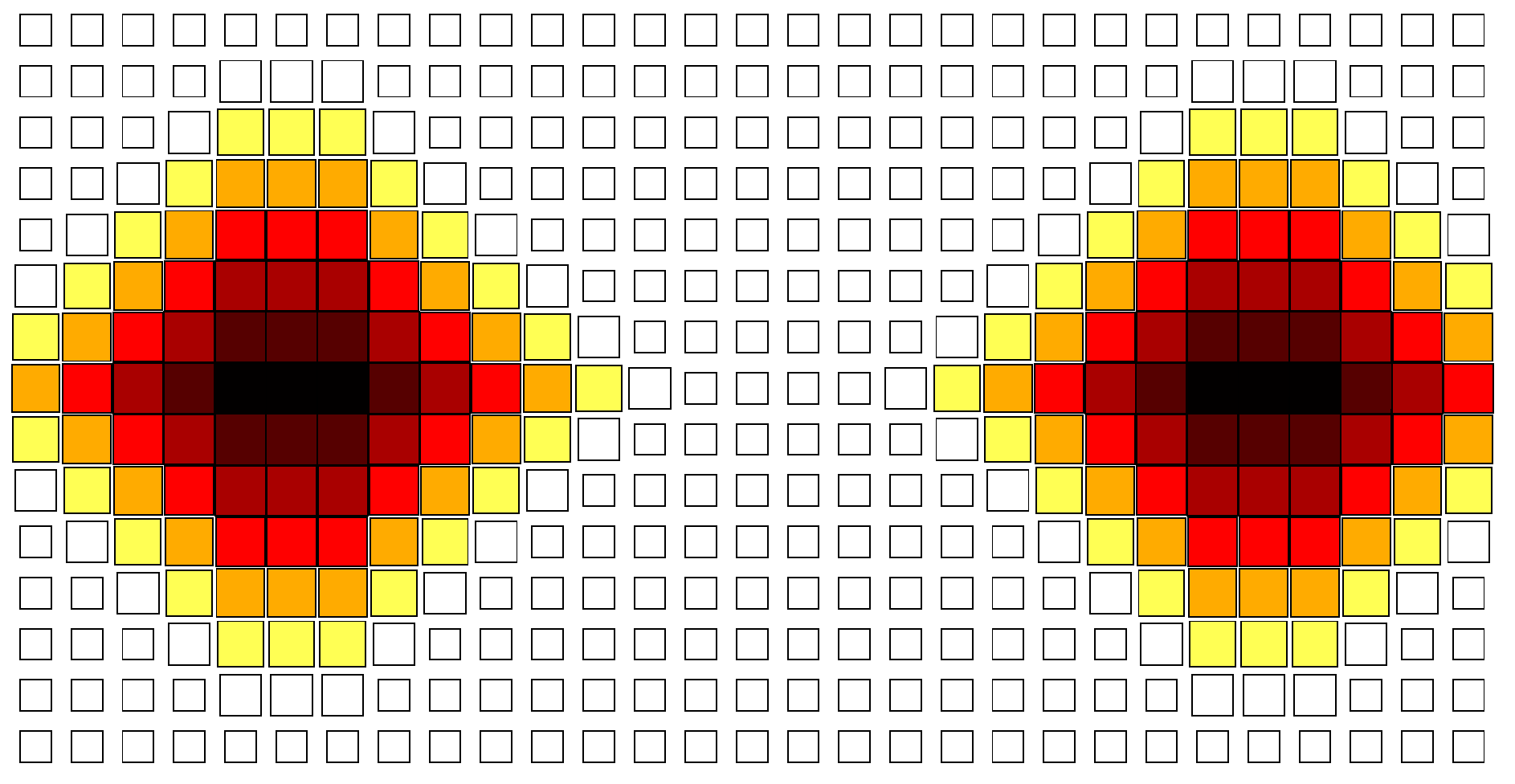}
} &
\subfigure[$N^{+}_i \ge \lfloor 2 N^{-}_i / 3 \rfloor$, $\Delta t_e = 10.$ ]{\label{fig:VortexPol_N23_10} 
  \includegraphics[width=0.45\textwidth]{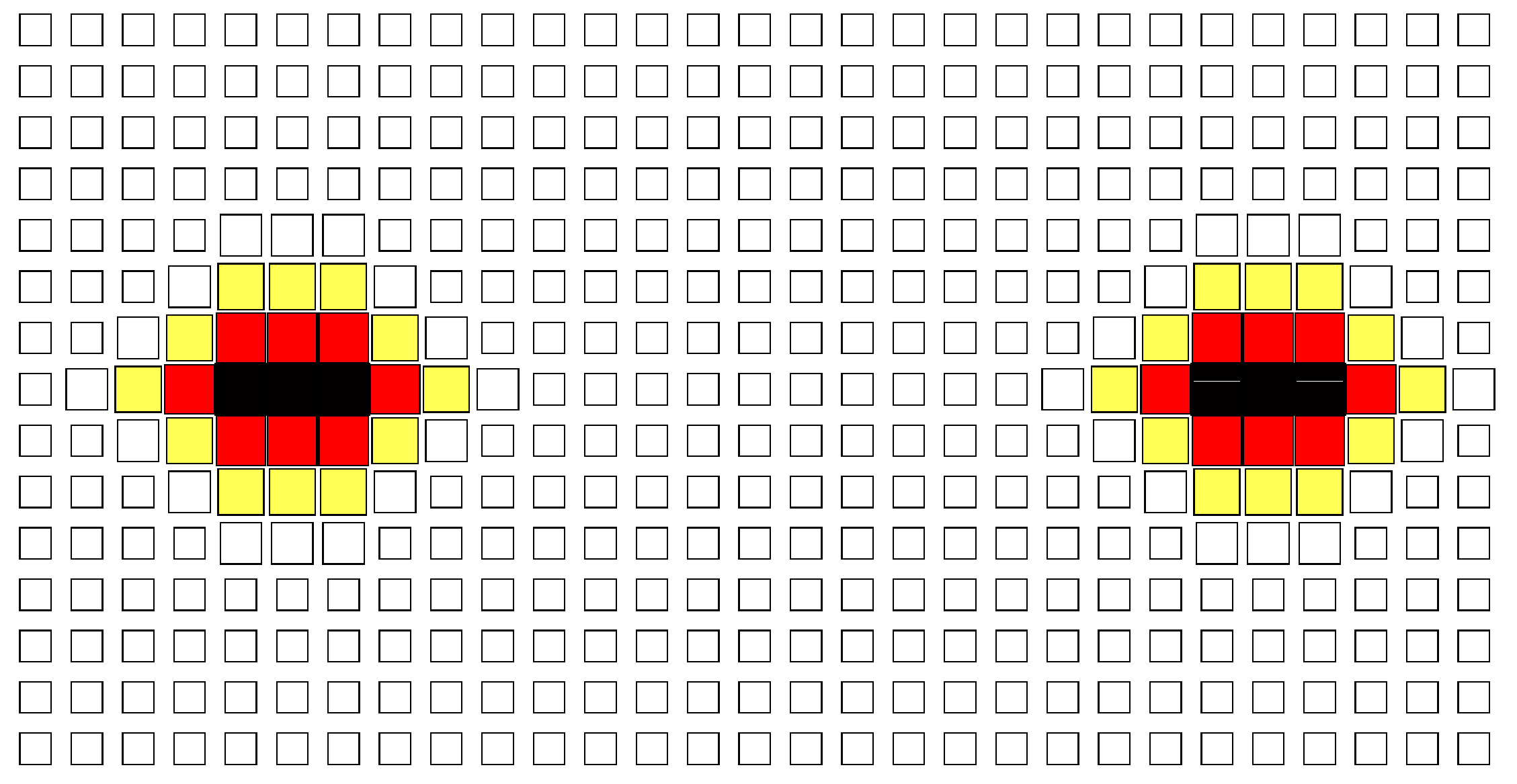}}
\\
\subfigure[$N^{+}_i \ge \lfloor N^{-}_i - 1   \rfloor$, $\Delta t_e = 5$.]{\label{fig:VortexPol_N_1_5} \includegraphics[width=0.45\textwidth]{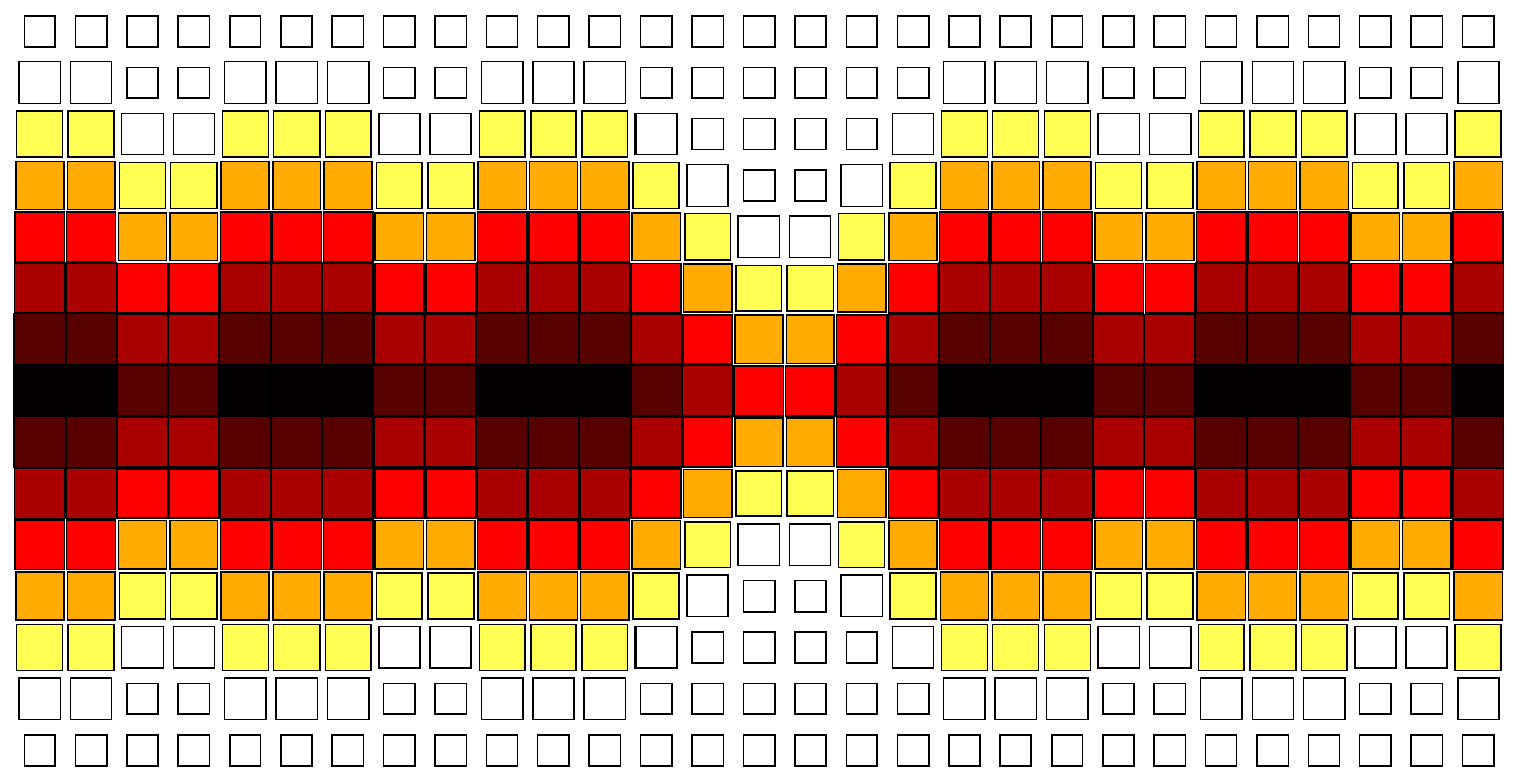}} &
\subfigure[$N^{+}_i \ge \lfloor 2 N^{-}_i / 3 \rfloor$, $\Delta t_e = 5.$ ]{\label{fig:VortexPol_N23_5} \includegraphics[width=0.45\textwidth]{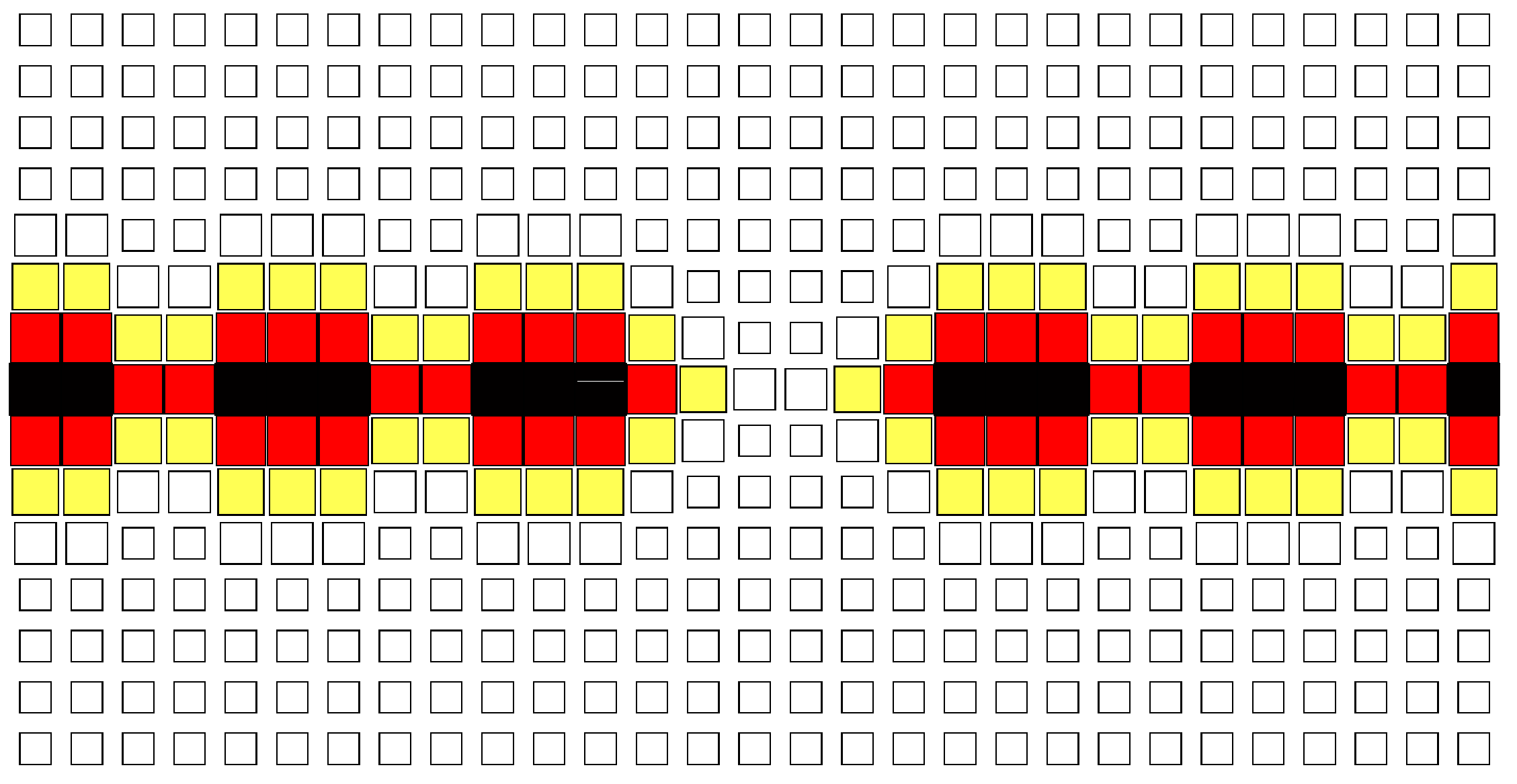}}
\\
\subfigure[$N^{+}_i \ge \lfloor N^{-}_i - 1   \rfloor$, $\Delta t_e = 1$.]{\label{fig:VortexPol_N_1_1} \includegraphics[width=0.45\textwidth]{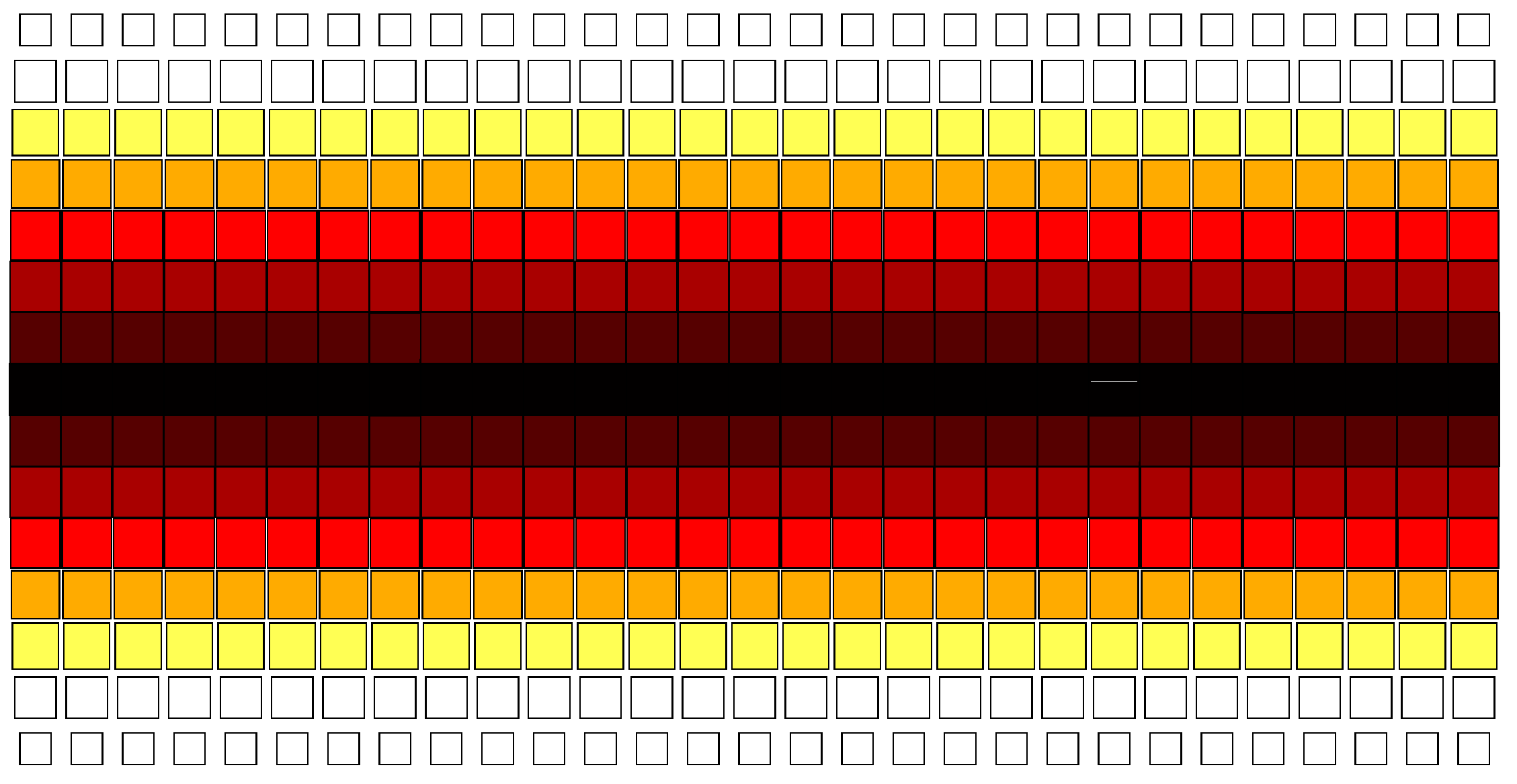}} &
\subfigure[$N^{+}_i \ge \lfloor 2 N^{-}_i / 3 \rfloor$, $\Delta t_e = 1.$ ]{\label{fig:VortexPol_N23_1} \includegraphics[width=0.45\textwidth]{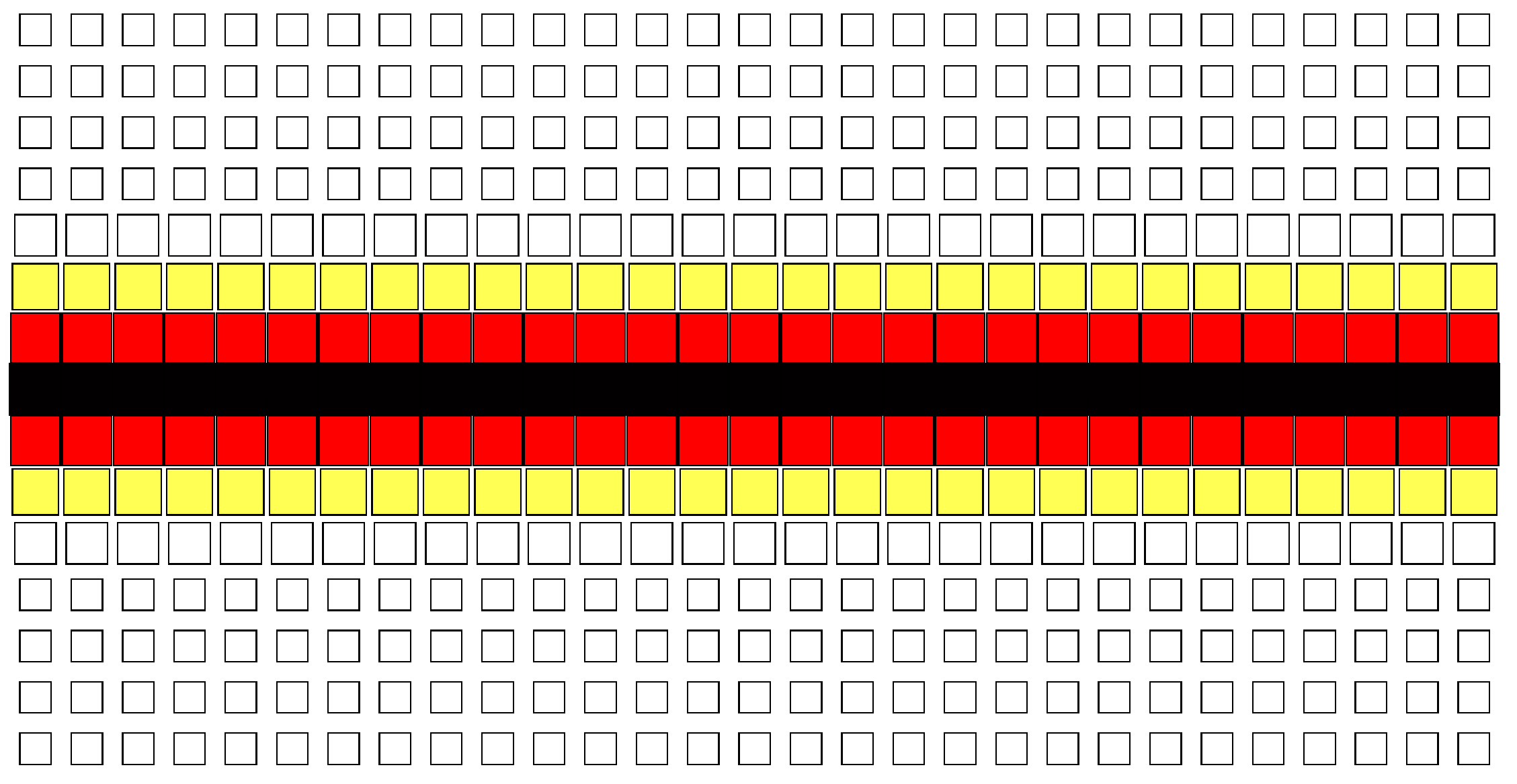}}

\end{tabular}

\caption{Average polynomial degree distribution for static $p$-adaptation with $\tau_{\max} = 10^{-1}$ and different estimation intervals.} \label{fig:VortexPol}
\end{figure}

The number of degrees of freedom of the $p$-adaptive simulations that obey the polynomial degree jump condition (b) \eqref{eq:Uadapt:N_1} is higher than for condition (a) \eqref{eq:Uadapt:N23}.
This is expected since many more elements are enriched in the former, as can be observed in Figure \ref{fig:VortexPol}.
The additional enrichment translates to computation times up to $25 \%$ higher when using condition (b) \eqref{eq:Uadapt:N_1}.

\begin{figure}[h]
\begin{center}

\includegraphics[width=0.55\textwidth]{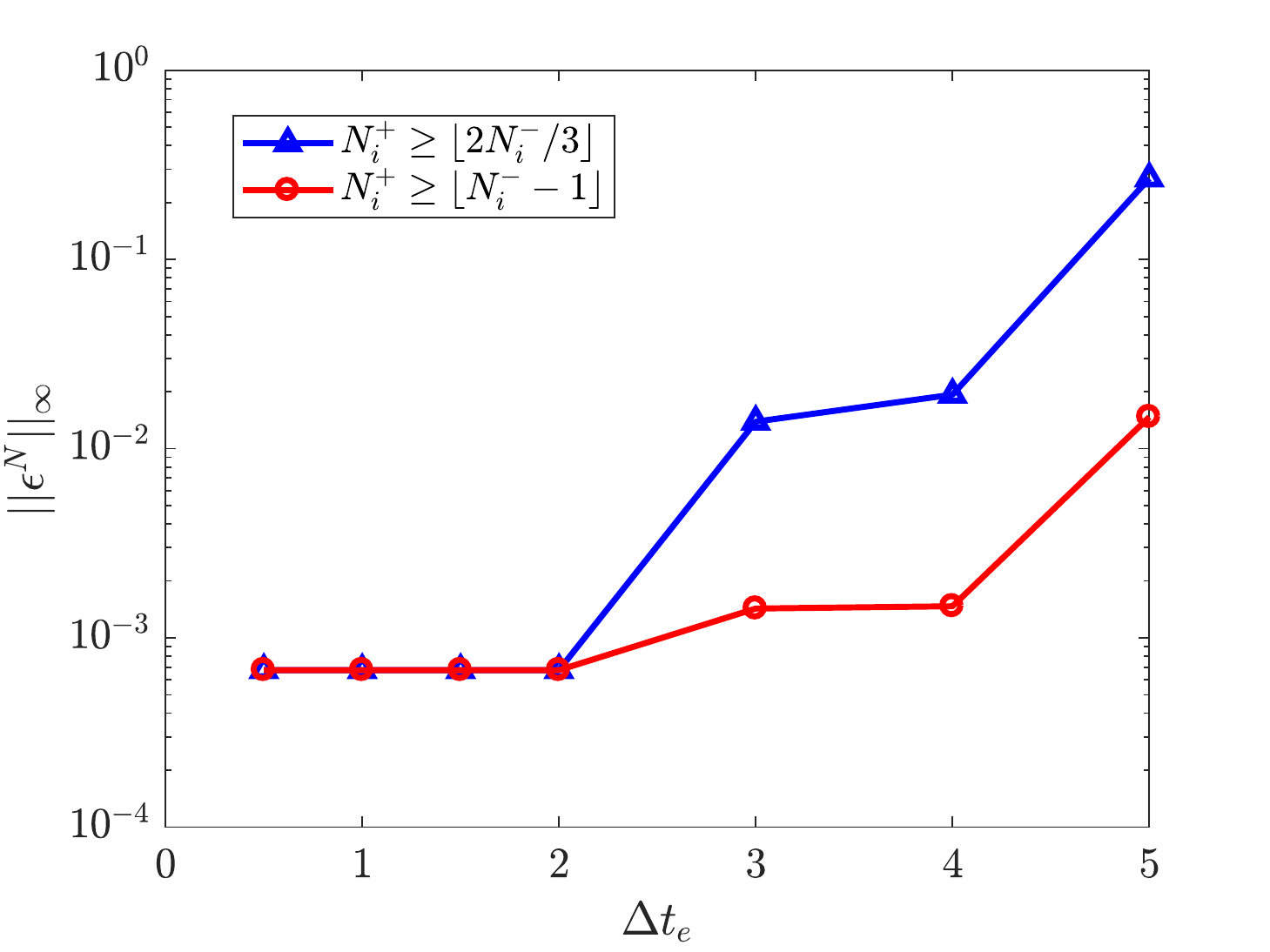} 
\caption{Combined effect of the polynomial degree jump condition and the estimation interval in the static $p$-adaptation error with threshold $\tau_{\max} = 10^{-1}$.}\label{fig:VortexOrderAc}
\end{center}
\end{figure}

An additional difference between the two polynomial degree jump conditions, which can be inferred from Figure \ref{fig:VortexPol}, is that condition (a) is more sensitive to the estimation/adaptation interval, $\Delta t_e$.
On the one hand, in the dynamically $p$-adaptive simulations and for a given $\Delta t_e$, it is more likely that the density pulse escapes the refined area for condition (a) than for (b), and arrives at an area where no $\tau$-estimation is possible ($P=1$) or where no extrapolation is possible ($P < 3$).
On the other hand, in statically $p$-adaptive simulations and for a given $\Delta t_e$, the refinement areas are more likely to be connected if condition (b) is used instead of (a).
This behavior is also illustrated in Figure \ref{fig:VortexOrderAc} for the static $p$-adaptation with the threshold $\tau_{\max} = 10^{-1}$.

\FloatBarrier

\subsection{Subsonic Flow Past a Cylinder}
\label{sec:Uadapt:flowpastcylinder}
We simulate the flow around a circular cylinder at a Reynolds number of $\Re = 100$ and a Mach number of $\Ma=0.15$ on a high-order curved ($P=3$) mesh with 1282 quadrilateral elements and the DGSEM method. 
We assess the performance of the truncation error-based static and dynamic $p$-adaptation methods and show that the static $p$-adaptation algorithm performs well in this example since the solution is statistically steady, as in most external aerodynamic problems.

\begin{figure}[h]
\begin{center}
\subfigure[Instantaneous horizontal velocity contours.]{\label{fig:Cyl_Re100_Contour} \includegraphics[width=0.46\textwidth]{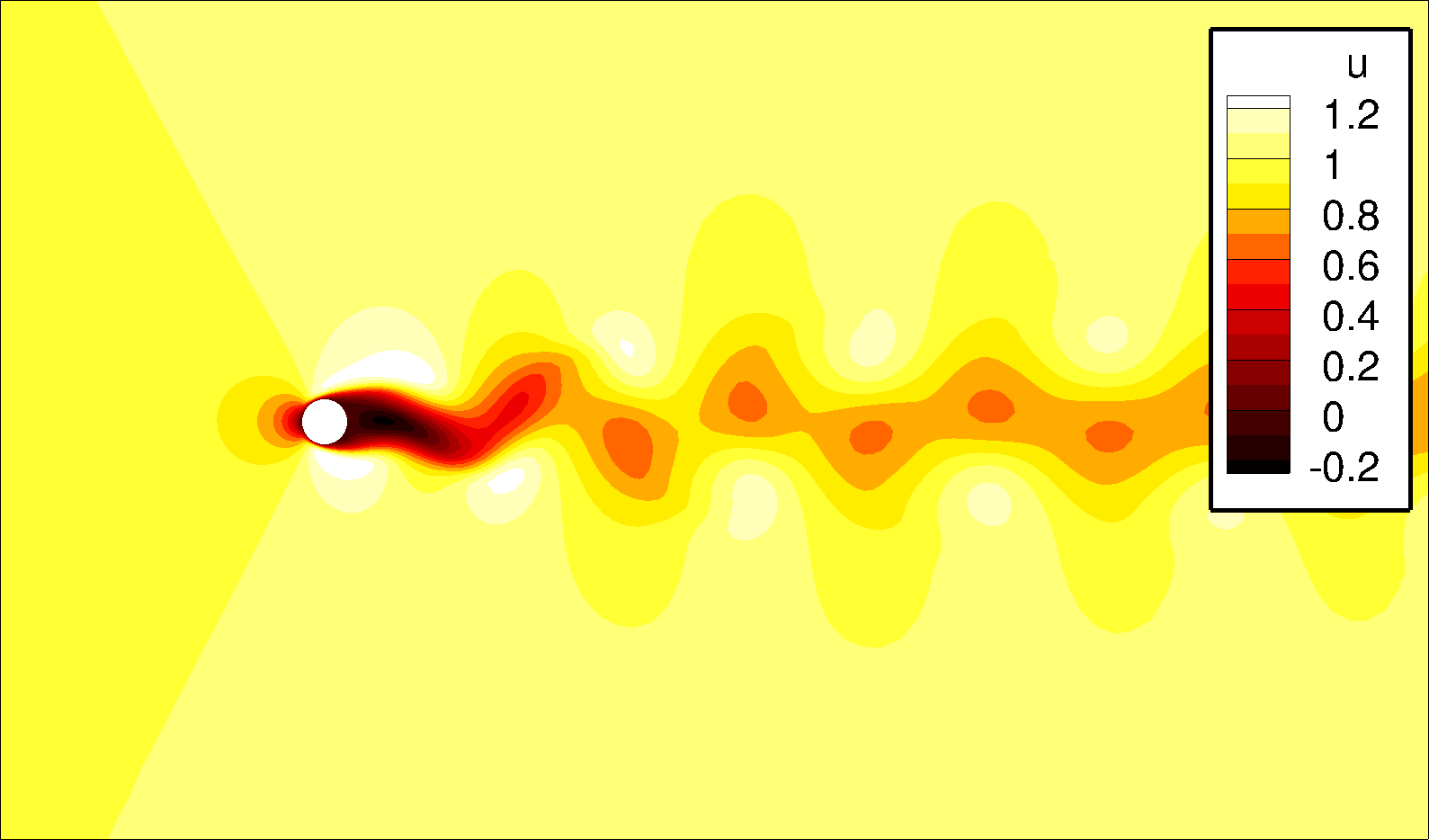}}\qquad
\subfigure[Instantaneous average polynomial degrees (dynamic $p$-adaptation) for $\tilde{\tau}_{\max} = 5$.]{\label{fig:Cyl_Re100_PolOrders} \includegraphics[width=0.46\textwidth]{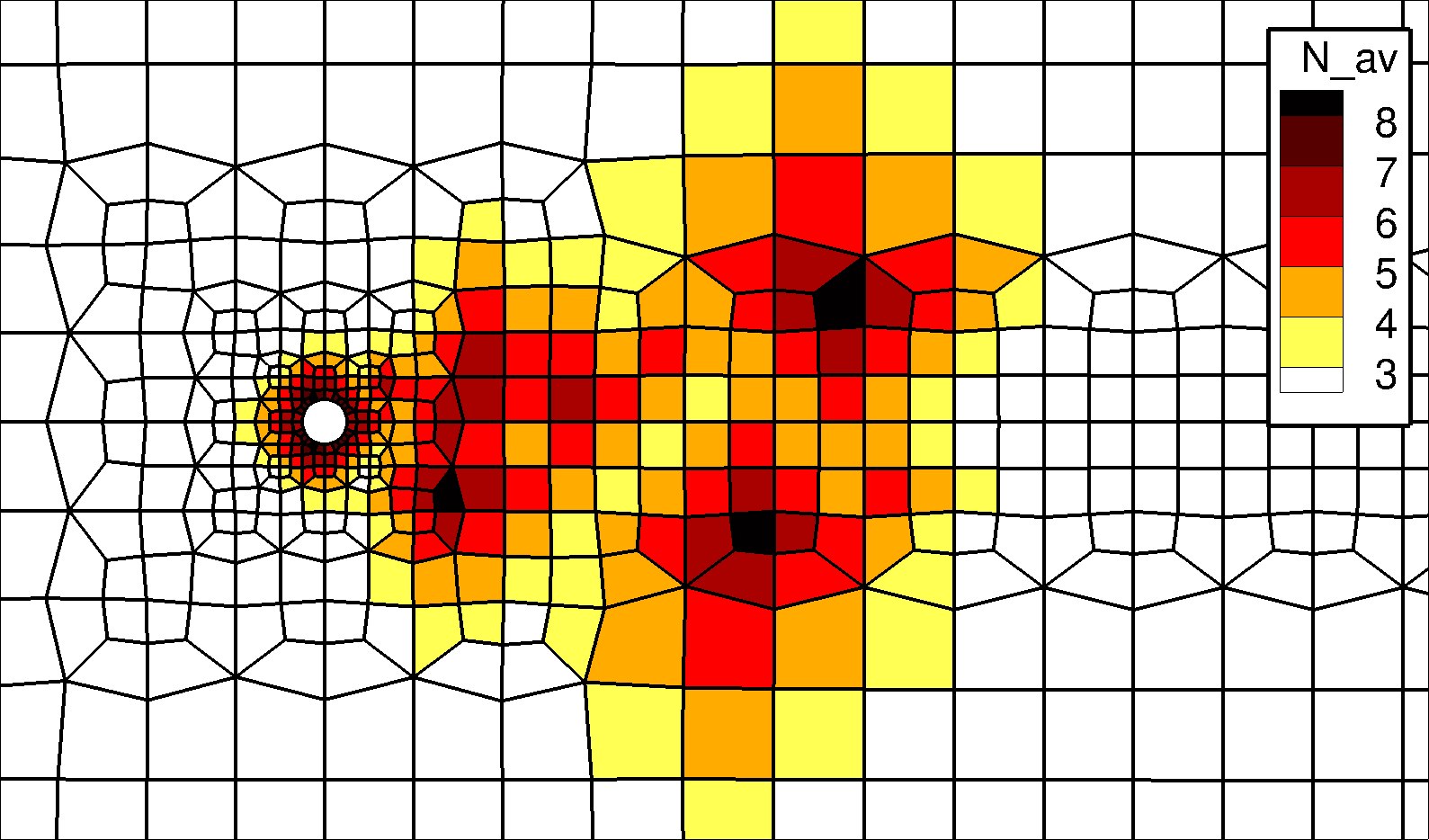}}

\caption{Vortex shedding behind a cylinder at $\Re=100$.} \label{fig:Cyl_Re100}
\end{center}
\end{figure}

Figure \ref{fig:Cyl_Re100} shows the mesh that was used, the instantaneous horizontal velocity contours, and an instantaneous distribution of polynomial degrees for the dynamic $p$-adaptation method.

The results presented in this section were obtained using a 40-core 2.10GHz Intel(R) Xeon(R) Gold 6230 CPU with 170 GB of RAM.
Each simulation was run with 10 cores and shared memory parallelization (OpenMP + guided schedule) to compute the spatial terms. 
Note that this parallel implementation has a near-optimal scalability for $p$-anisotropic discretizations and the selected OpenMP schedule, as was shown in \cite{RuedaRamirez2019a}.
We remark that the guided OMP schedule acts directly as a dynamic load balancing technique in the simulations with dynamic $p$-adaptation.

\begin{figure}[h]
\begin{tabular}{cc}

\multicolumn{2}{c}{
	\includegraphics[width=0.45\textwidth]{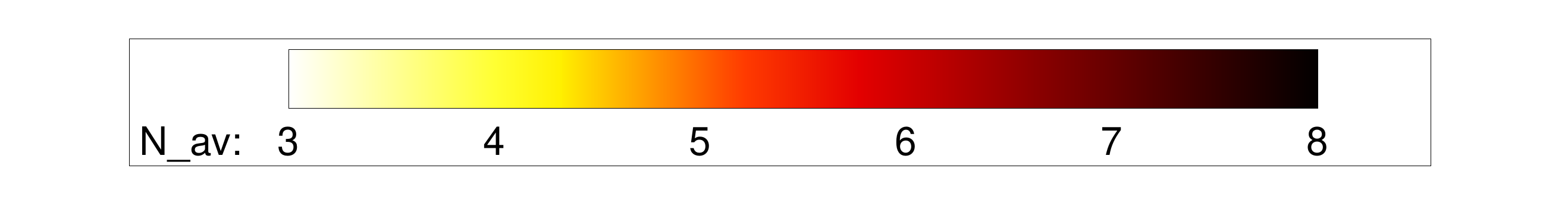}
} \\

\subfigure[Trad. $\tau_{\max} = 10^{-1} (\NDOF = 43574)$.]
{
  \label{fig:CylPol_N_1_10} 
  \includegraphics[width=0.45\textwidth]{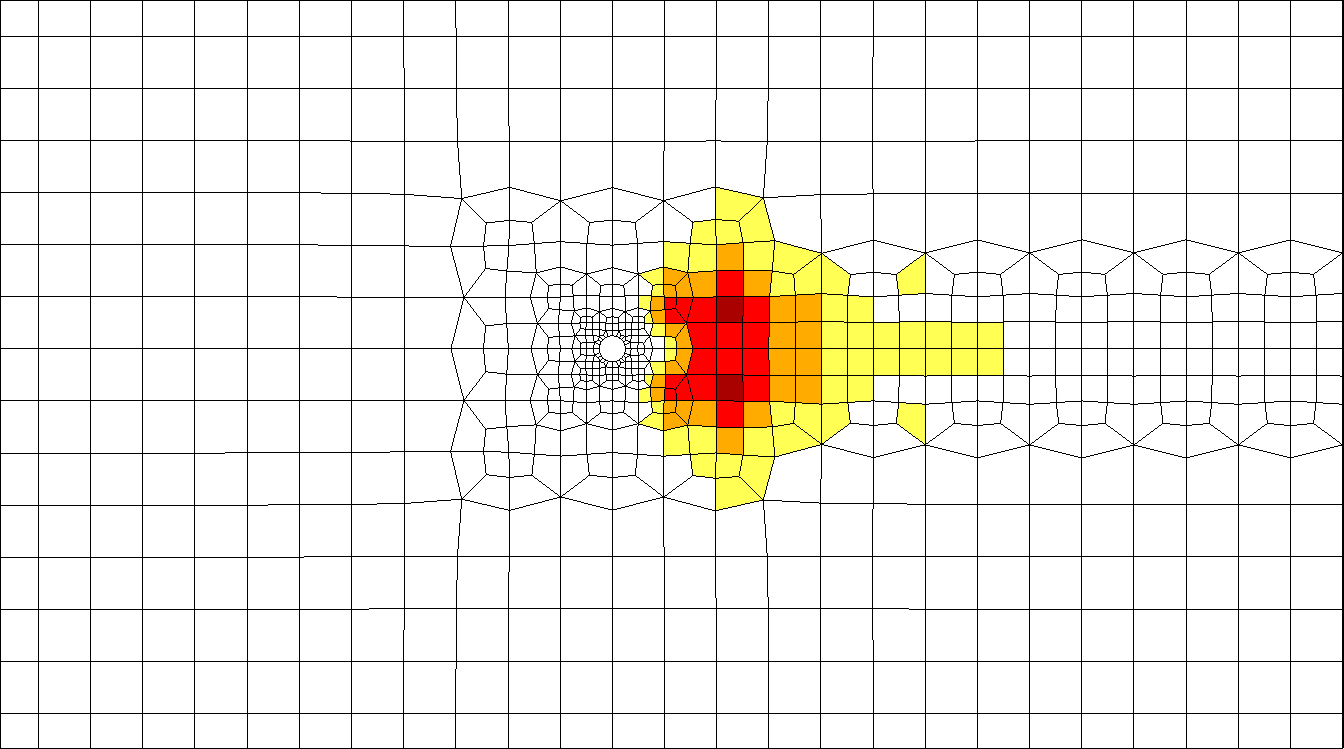}
} &
\subfigure[New. $\tilde{\tau}_{\max} = 10^{2} (\NDOF = 41808)$.]{\label{fig:CylPol_N23_10} 
  \includegraphics[width=0.45\textwidth]{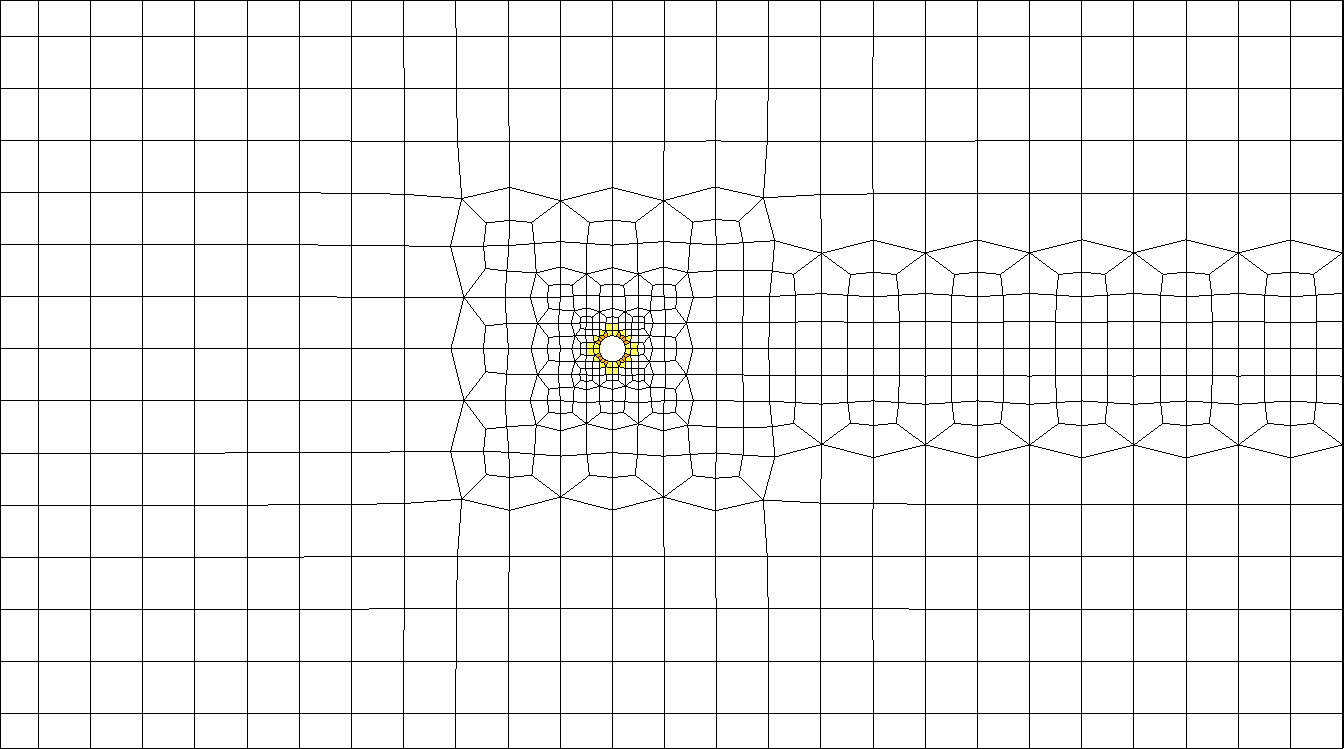}}
\\
\subfigure[Trad. $\tau_{\max} = 10^{-2} (\NDOF = 69644)$.]{\label{fig:CylPol_N_1_5} \includegraphics[width=0.45\textwidth]{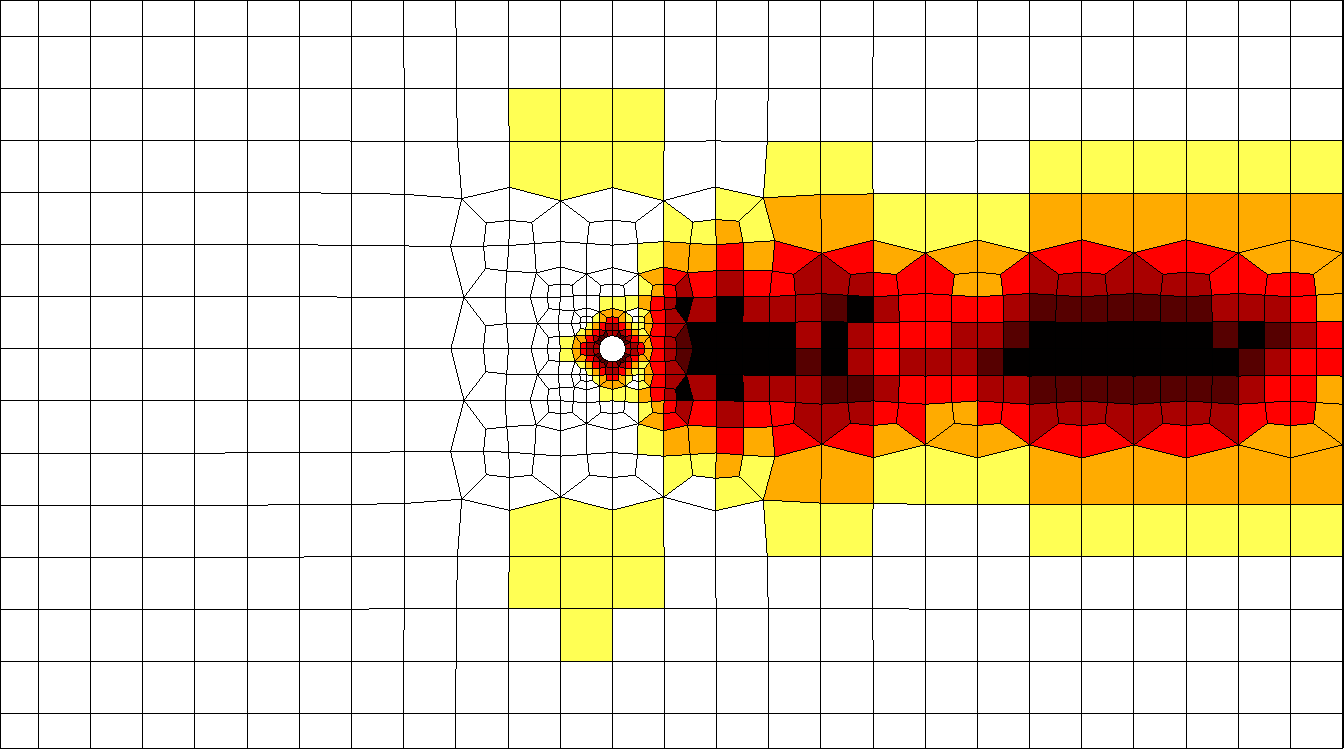}} &
\subfigure[New. $\tilde{\tau}_{\max} = 1 (\NDOF = 67766)$.]{\label{fig:CylPol_N23_5} \includegraphics[width=0.45\textwidth]{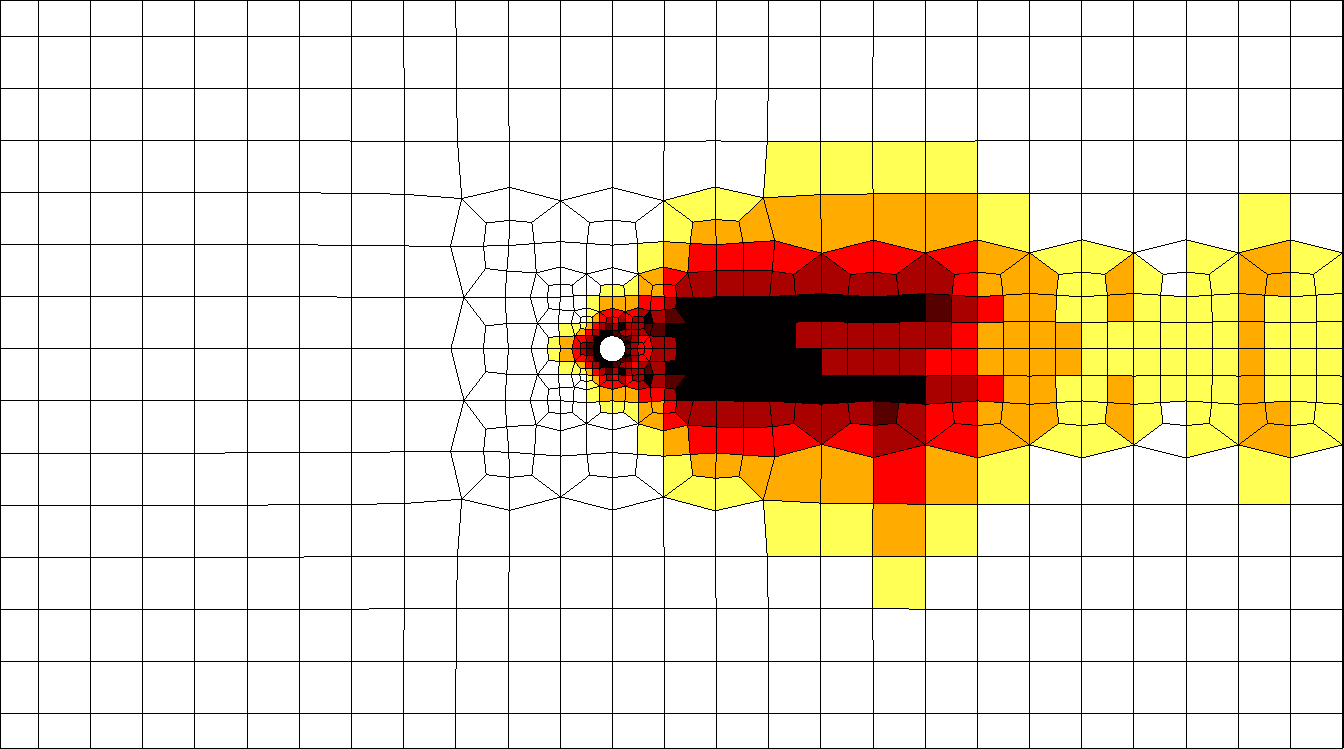}}
\\
\subfigure[Trad. $\tau_{\max} = 5 \times 10^{-3} (\NDOF = 90192)$.]{\label{fig:CylPol_N_1_1} \includegraphics[width=0.45\textwidth]{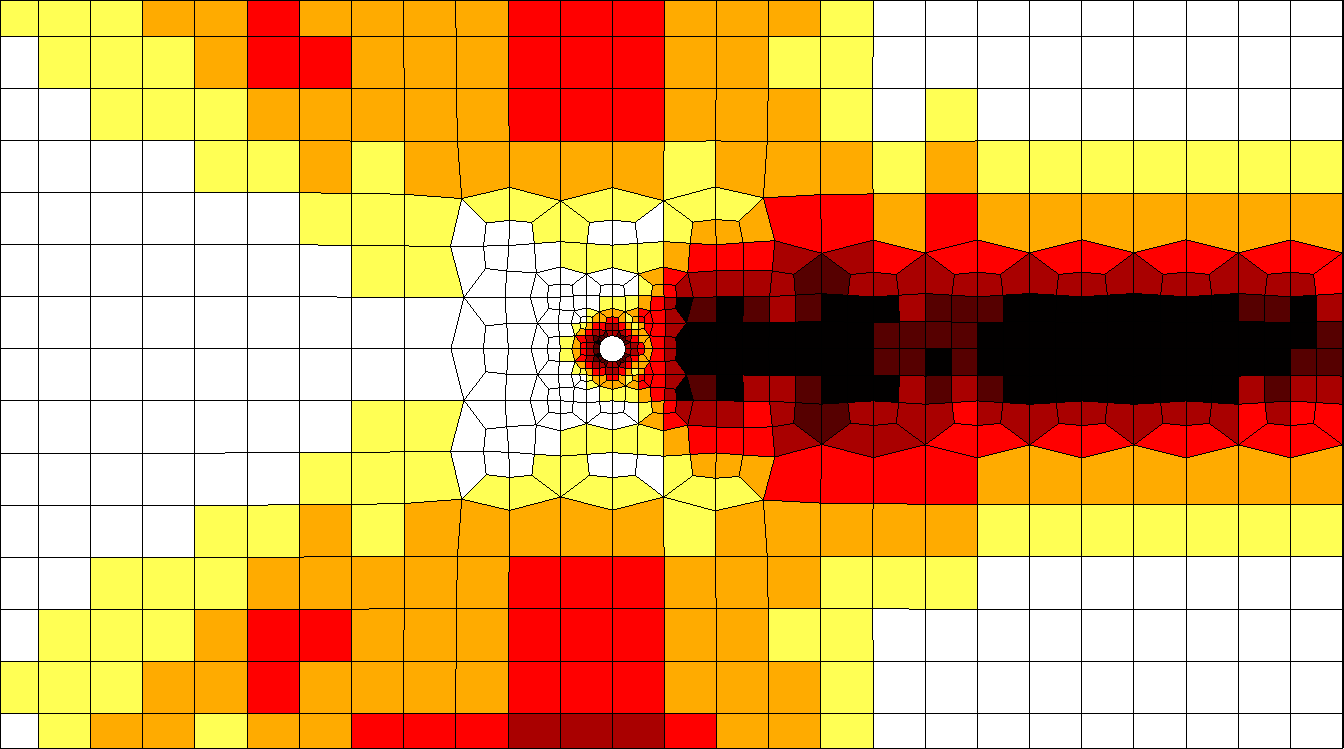}} &
\subfigure[New. $\tilde{\tau}_{\max} = 10^{-1} (\NDOF = 95930)$.]{\label{fig:CylPol_N23_1} \includegraphics[width=0.45\textwidth]{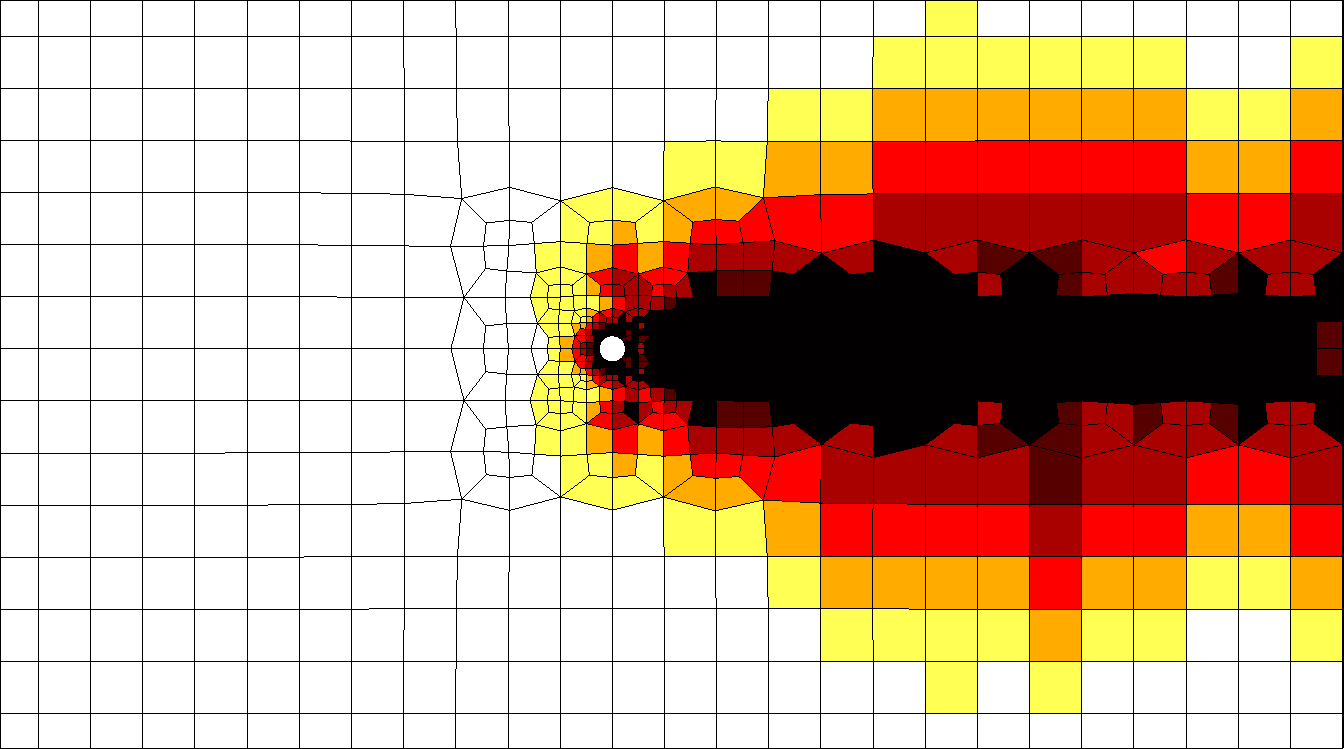}}

\end{tabular}

\caption{Comparison of the two possible formulations of the truncation error: traditional (left) and new (right). Average polynomial degree distribution for static $p$-adaptation for different error thresholds.} \label{fig:CylPol}
\end{figure}

For the $p$-adaptive simulations, the new form of the truncation error is retained because it was shown to work more efficiently on the lift and drag predictions than the traditional form.
The reason why is easily seen in Figure \ref{fig:CylPol}, which shows the contours of the average polynomial degrees for both formulations of the truncation error as the error threshold, $\tau_{\max}$, is reduced in a static $p$-adaptation method.
For a similar number of degrees of freedom, the $p$-adaptation algorithm that uses the traditional form of $\tau$ tends to enrich large elements that are away from the cylinder, whereas the new form tends to enrich only the boundary layer area and the wake.
As explained in Section \ref{sec:Uadapt:NewTruncError}, the main difference between the two approaches is the weight they assign to the volume of each element.

The truncation error-based static and dynamic $p$-adaptation algorithms are tested with truncation error thresholds ranging between $10^{-1} \le \tilde{\tau}_{\max} \le 10^2$, and four estimation/adaptation intervals $\Delta t_e = 0.5, 1, 3, 6$ in non-dimensional time units, taking into account that the vortex shedding period is expected to be $T=6$.
Furthermore, the polynomial degree jump across faces is limited to $N^{+}_i \ge \lfloor N^{-}_i - 1 \rfloor$ \eqref{eq:Uadapt:N_1}, since this condition provides robustness to the simulation and allows larger estimation intervals, as discussed in Section \ref{sec:FlowVortex}.
Additionally, the maximum polynomial degree was set to $N_{\max}=8$ and the minimum polynomial degree to $N_{\min} = 3$.
This minimum polynomial degree allows the dynamic $p$-adaptation to always have enough points to perform the directional truncation error extrapolation.

In the dynamic $p$-adaptation algorithm, the sub-meshes that are used for the truncation error estimation are constructed every $\Delta t_e$ time units.
After that, the error is estimated using the $\tau$-estimation method and the polynomial degrees are changed accordingly.
At every adaptation stage, we only allow the polynomial degree to decrease by one in each element to reduce spurious oscillations that may arise because of large polynomial degree jumps throughout the simulation (a phenomenon discussed in Section \ref{sec:DynamicPAdapt}).
In the static $p$-adaptation algorithm, the $\tau$-estimation sub-meshes are only constructed once at the beginning of the simulation. 
Thereafter, an estimation simulation with polynomial degree $P=4$ is run for a sampling time of $T_e = 12$, i.e. two vortex shedding cycles.
The polynomial degrees are then adapted using \textit{strategy 2} \eqref{eq:Uadapt:totalTruncMap}, and the rest of the simulation runs without further modifications.

\begin{figure}[h]
\begin{center}
\begin{tabular}{cc}

\multicolumn{2}{c}{
	\includegraphics[width=0.35\textwidth]{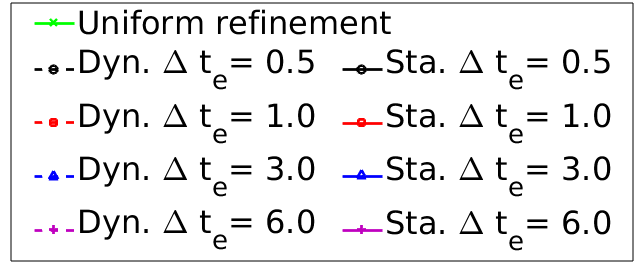}
} \\

\subfigure[Absolute lift vs. NDOF.]{\label{fig:Cyl_Lift_DOFs} \includegraphics[width=0.46\textwidth]{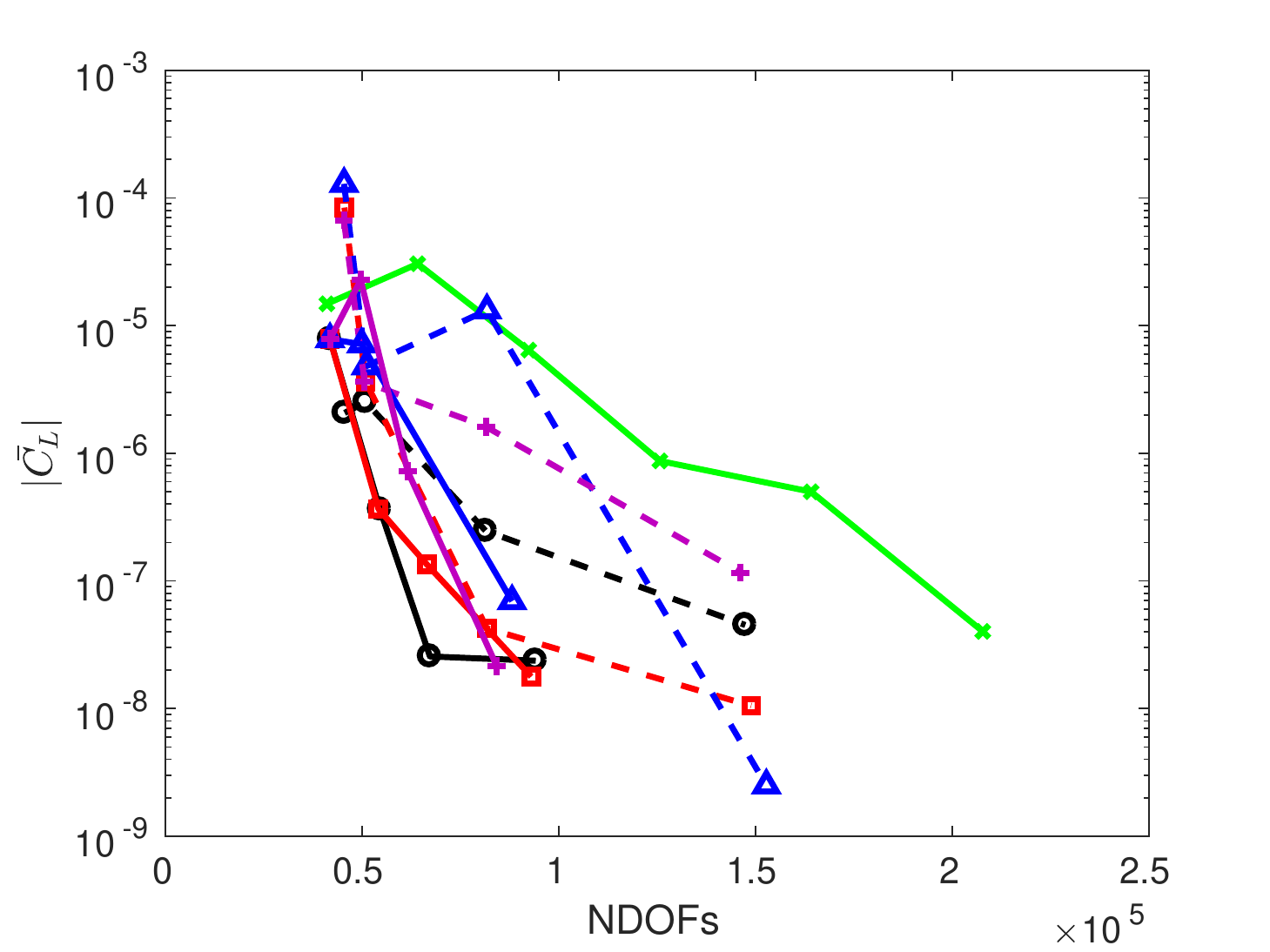}}
&
\subfigure[Absolute lift vs. CPU-Time.]{\label{fig:Cyl_Lift_Time} \includegraphics[width=0.46\textwidth]{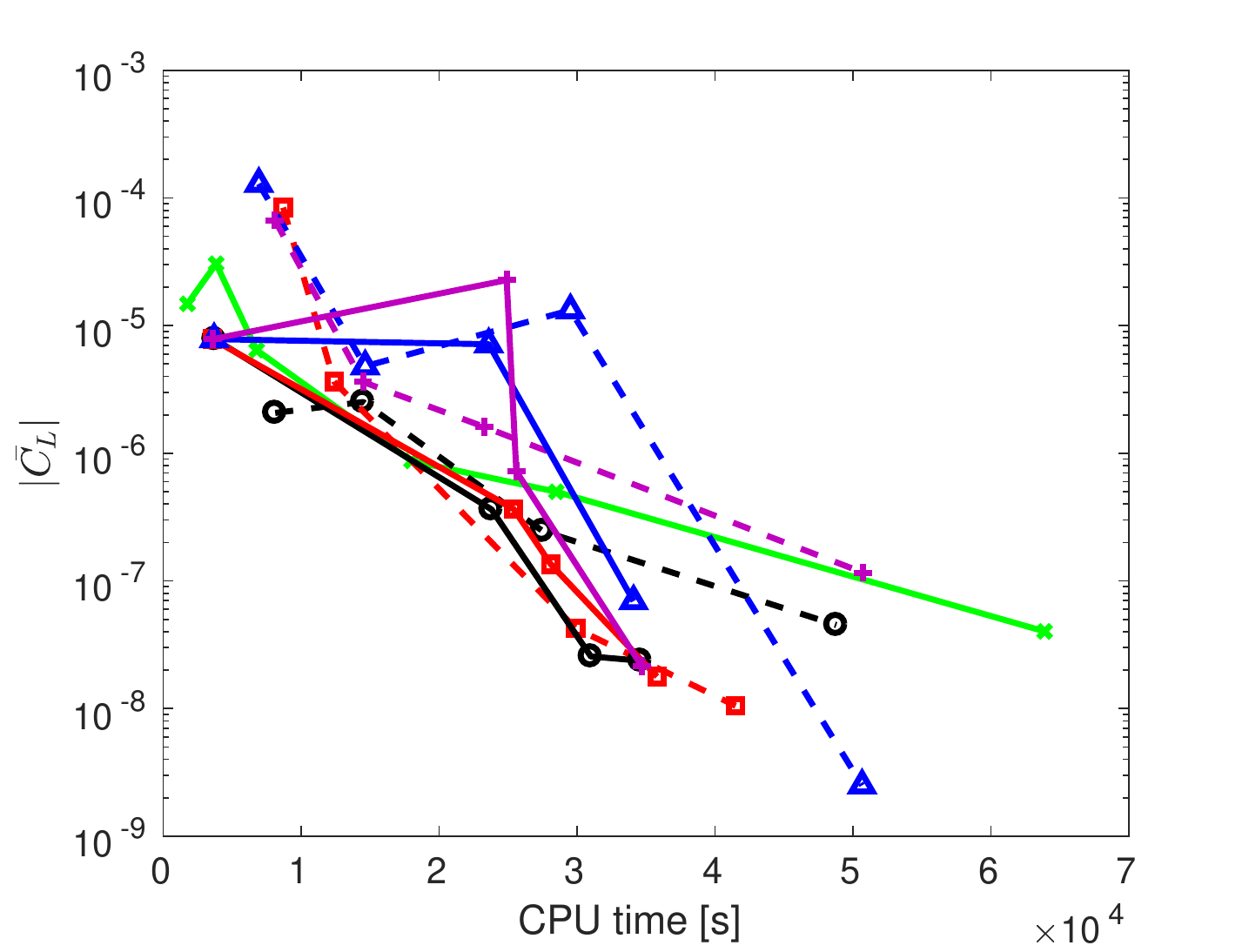}}
\\

\subfigure[Drag error vs. NDOF.]{\label{fig:Cyl_Drag_DOFs} \includegraphics[width=0.46\textwidth]{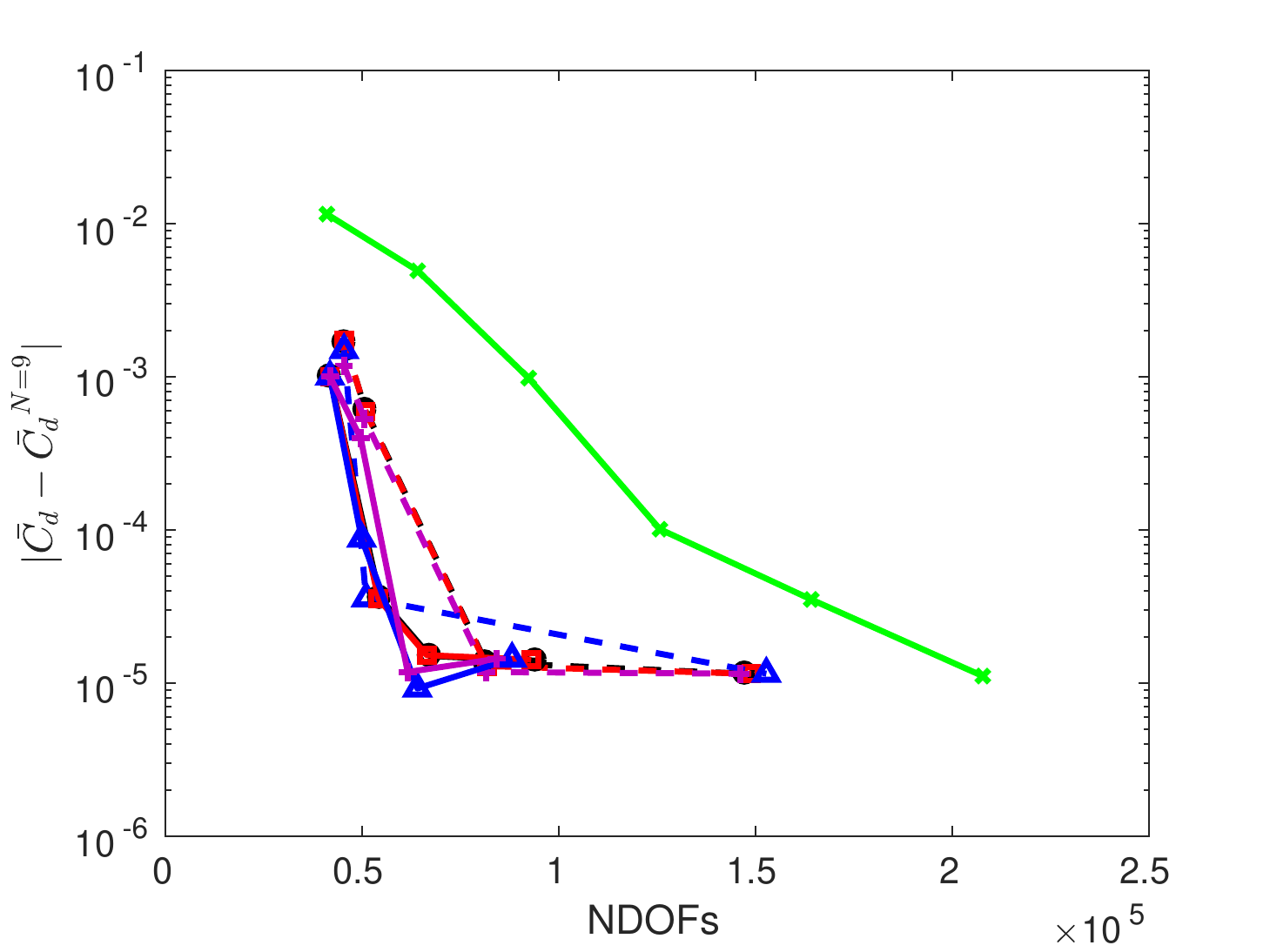}}
&
\subfigure[Drag error vs. CPU-Time.]{\label{fig:Cyl_Drag_Time} \includegraphics[width=0.46\textwidth]{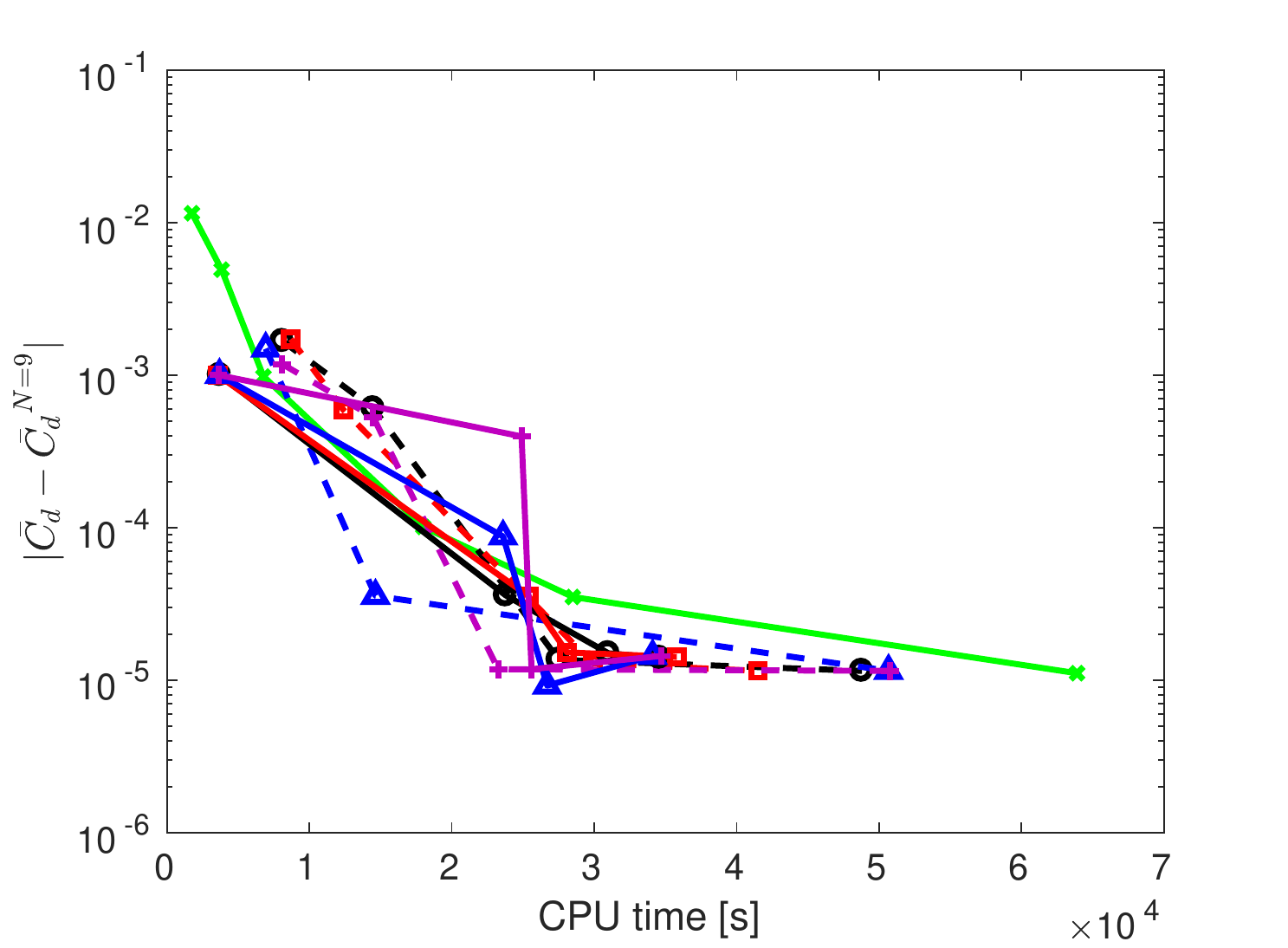}}

\end{tabular}

\caption{Performance of the static and dynamic $p$-adaptation procedures for the flow past a cylinder at $\Re = 100$, $\Ma = 0.15$.
Each point in the graph corresponds to a truncation error threshold, ranging between $10^{-1} \le \tilde{\tau}_{\max} \le 10^2$. 
Different intervals between adaptation/estimation stages are represented with lines of different colors (black, red, blue and purple).
Performance of the uniform polynomial degree refinement (in green) is also included.
} \label{fig:Cyl_Performance}
\end{center}
\end{figure}

Figure \ref{fig:Cyl_Performance} shows the performance of the uniform $p$-refinement, dynamic and static truncation error-based $p$-adaptation algorithms.
The mean absolute lift and the mean drag error 
(the latter with respect to a solution of order $N=9$) 
are plotted as a function of the number of degrees of freedom (NDOF) and the computation time for each of the methods.
Lift and drag are monitored for 100 time units and their average values are computed. 
The reported computation time is the sum of the CPU time that is needed to advance 100 time units and the CPU time that is needed for the estimation.

Because of the adaptive time-stepping, the monitored variables had to be re-sampled at a uniform time-step sequence to calculate the averages.
A brief description of the process to obtain the re-sampled data is provided in \ref{app:Uadapt:postpro}.
In addition, the number of degrees of freedom that is shown for the dynamic $p$-adaptation simulations corresponds to a weighted average,
\begin{equation}
\NDOF_{\mathrm{dyn}} = \frac{1}{S} \sum_{i=1}^S \NDOF_i,
\end{equation}
where $S$ is the number of simulation time steps and $\NDOF_i$ corresponds to the number of degrees of freedom of the discretization in the iteration $i$.\\

It can be observed in Figure \ref{fig:Cyl_Performance} that the truncation error-based locally adaptive simulations need fewer degrees of freedom than the simulations with uniform order.
Furthermore, in contrast to the advected pulse example, the static $p$-adaptation method needs fewer degrees of freedom than the dynamic $p$-adaptation method for the same accuracy.
The main reason for this behavior is that the dynamic $p$-adaptation algorithm might overestimate the polynomial degree needed when the polynomial degree of the reference mesh (used for the estimation), $P$, is low (a behavior discussed in Section \ref{sec:DynamicPAdapt}).
In fact, the dynamic $p$-adaptation algorithm is more likely to overpredict the required polynomial degree than the static algorithm since the minimum specified polynomial degree acts sometimes as the estimation polynomial degree in dynamically $p$-adaptive simulations, $P=N_{\min}=3$, which is lower than the estimation polynomial degree of the static $p$-adaptation algorithm, $P=4$.

The behavior of the error with respect to the computation times is highly dependent on the implementation, the hardware used, and the problem.
The results obtained with the current implementation in HORSES3D \cite{ferrer2022horses3d} are reported as a reference.
As can be observed, the performance is different for each variable analyzed, but in general a speed-up of about $2.2$ can be observed for the static $p$-adaptation algorithm at the highest level of accuracy that is reached.
The dynamic $p$-adaptation algorithm has the same performance as the static $p$-adaptation algorithm in some cases, and in some others it exhibits a worse performance.
The main reason for that is that the dynamic $p$-adaptation algorithm is more sensitive to the estimation interval, $\Delta t_e$, and that it may also suffer from non-physical oscillations in the solution and its gradients due to the frequent jumps in the polynomial degree.

The truncation error-based $p$-adaptation methods show the best CPU-time performance when measuring the mean absolute lift, where speed-ups can be observed in virtually all the error range considered for small enough $\Delta t_e$.
When measuring the mean drag error, the truncation error-based $p$-adaptation performs relatively similar to the uniform $p$-refinement with respecto to CPU-time (if $\Delta t_e$ is small enough) down to an error of $|\bar{C}_d - \bar{C}_d^{N=9}| \approx 5 \times 10^{-5}$.
Below that error, the truncation error-based $p$-adaptation algorithms outperform the uniform refinement technique.

\FloatBarrier

\section{Conclusions} \label{sec:Uadapt:Conclustions}

In this paper, we have extended the truncation error-based $p$-adaptation method to unsteady problems. 
First, we presented a new form of the truncation error, which holds close similarities to the formulation traditionally used in the literature for variational methods \cite{Kompenhans2016,RuedaRamirez2019,RuedaRamirez2019a,RuedaRamirez2019b}. 
The new form of the truncation error performs well and similarly to the traditional form for a test case with uniform mesh size. However, when considering a nonuniform mesh size, the new form outperforms the traditional formulation. 
Second, we extended the $\tau$-estimation method to estimate the truncation error of unsteady flow problems with the DGSEM. 
The method developed here retains the anisotropic properties and the ability to be estimated in a multigrid cycle, as proposed by the authors in \cite{RuedaRamirez2019a}.
Third, we proposed two truncation error-based $p$-adaptation strategies: the dynamic and static adaptation methods. 
We analyzed both strategies and used them successfully to enhance the performance of DGSEM in the open-source framework HORSES3D \cite{ferrer2022horses3d}. 
We conclude that the static $p$-adaptation method performs better than the dynamic one in statistically steady problems where the flow features are concentrated in a small part of the domain. Similarly, the dynamic $p$-adaptation method outperforms the static one when the flow features move through a large portion of the domain. 
For the test cases considered here, significant speed-ups (up to $\times 4$) are reported.

\appendix
\section{A Note on Post-Processing}
\label{app:Uadapt:postpro}
In this section, we provide a short description of the post-processing method used to acquire the results presented in Section \ref{sec:Uadapt:flowpastcylinder}. 
Since we carried out the calculations with a constant CFL (instead of a constant time-step size) and stored the lift and drag at every time step, we re-sample the lift and drag signals to obtain equispaced data in time.
The process consists of three steps:

\begin{enumerate}
    \item We take the last part of the signal to avoid the effect of any transients from the restart. 
    For all simulations considered, the last $20$ time units (of a total of $100$) showed to have a periodic behavior. 
    
    \item We cut the signal from left and right to ensure that we are averaging over entire periods of the signal.
    First, we compute the mean value of the signal (lift or drag), locate the first position where this mean value appears in the time series, and remove data left from that point. 
    Then, we locate the last position where this mean value appears in the time series with a slope of the same sign and remove data right from that point.
    The resulting signal has $n$ points.
    
    \item We feed the time series obtained in step 2 into the MATLAB function \textit{resample} to get a new signal with $4 \times n$ equidistant points. 
    
\end{enumerate}
\section*{Acknowledgments}

AR acknowledges funding through the Klaus-Tschira Stiftung via the project "HiFiLab".
EF and GN acknowledge the financial support of the European Union’s Horizon 2020 research and innovation programme under the Marie Skłodowska-Curie grant agreement (MSCA ITN-EID-GA ASIMIA No 813605). 
GR and EV acknowledge the funding received by the Grant SIMOPAIR (Project No. RTI2018-097075-B-I00) funded by MCIN/AEI/ 10.13039/501100011033 and by ERDF A way of making Europe. 
AR, EV and EF thank the European Union Horizon 2020 Research and Innovation Program under the Marie Sklodowska-Curie grant agreement No 675008 for the SSeMID project.
Finally, all authors gratefully acknowledge the Universidad Politécnica de Madrid (www.upm.es) for providing computing resources on the Magerit Supercomputer.

\bibliography{Biblio.bib}

\end{document}